\documentclass[12pt]{amsart}
 
\usepackage{amsfonts, 
                    amsmath, 
                    amssymb, 
                    mathrsfs,
                    mytitle,
                    amscd, 
                    verbatim,
                   graphicx }

\usepackage[T1]{fontenc} 

\input{epsf.sty}

\DeclareMathAlphabet{\mathpzc}{OT1}{pzc}{m}{n}

\newcommand{\mathsmall}[2]{\mbox{\fontsize{#1}{1}$#2$}}

\def\Rqs{{\mathbb L}}
\def\fRqs{{\overline{\mathbb L}}}
\def\RqsFN{{\mathbb L}} 

\def\Am{\mbox{\boldmath$A$}} 
\def\Bm{\mbox{\boldmath$B$}} 
\def\Vm{\mbox{\boldmath$V$}} 
\def\Im{\mbox{\boldmath$\mathrm I$}} 
\def\Imo{\mbox{\boldmath$\overline{\mathrm I}$}} 
\def\Vmh{\mbox{\boldmath$V_{\hbar}$}} 
\def\Ql{\mbox{\boldmath$Q$}} 
\def\W{\mbox{\boldmath$W$}} 

\def\sAm{\mbox{\small \boldmath$A$}} 
\def\sBm{\mbox{\small \boldmath$B$}}

\def\sW{\mbox{\small \boldmath$W$}} 

\def\arh{\rule{0mm}{5mm}}

\newcommand\lkb[2]{\mbox{\boldmath$H$}_{#2,#1}}
\newcommand\lkbt[2]{\mbox{\boldmath$\widetilde H$}_{#2,#1}}

\def\uqsl{U_q({\mathfrak s}{\mathfrak l}_2)}
\def\sl{{\mathfrak s}{\mathfrak l}_2}

\def\Alg{\mbox{\large$\mathpzc U$}}

\def\invo{\mbox{\boldmath$\iota$}} 

\def\TLsi{\mbox{\textsf{\textit{TL}}}}
\def\Isi{\mbox{\textsf{\textit{I}}}}
\def\Jsi{\mbox{\textsf{\textit{J}}}}

\def\Rm{\mathsf R}

\def\Bas{{\mathcal W}}
\def\ABas{{\mathcal A}}
\def\BBas{{\mathcal B}}

\def\Bid{{\mathchoice {\rm {1\mskip-4.5mu l}} {\rm
{1\mskip-4.5mu l}} {\rm {1\mskip-3.8mu l}} {\rm {1\mskip-4.3mu l}}}}

\def\beq {\begin{equation}}
\def\eeq {\end{equation}}

\newtheorem{theorem}{Theorem}
\newtheorem{propos}[theorem]{Proposition}

\newtheorem{conjecture}[theorem]{Conjecture}

\newtheorem{defi}[theorem]{Definition}
\newtheorem{qes}[theorem]{Question}
\newtheorem{exe}[theorem]{Example}
\newtheorem{lemma}[theorem]{Lemma}

\newcommand{\chapter}[1]{{\bf\Large #1}}

\newcommand{\txtqt}[2]
{\mbox{\raise .15em\hbox{$#1$}\!\big/\!\raise -.18em\hbox{$#2$}}}

\def\rightarr{\mbox{\begin{picture}(15,1)(0,0)\put(0,0)
{\vector(1,0){15}}\end{picture}}}

\def\underrightarrow#1{\mathop
{\vtop{\ialign{##\crcr
$\hfil\displaystyle{#1}\hfil$\crcr\noalign{\kern1pt\nointerlineskip}
\rightarr\crcr\noalign{\kern.5pt}}}}\limits}

\def\leftarr{\mbox{\begin{picture}(15,1)(0,0)\put(15,0)
{\vector(-1,0){15}}\end{picture}}}

\def\underleftarrow#1{\mathop
{\vtop{\ialign{##\crcr
$\hfil\displaystyle{#1}\hfil$\crcr\noalign{\kern1pt\nointerlineskip}
\leftarr\crcr\noalign{\kern.5pt}}}}\limits}

\def\blm{\begin{lemma}}
\def\elm{\end{lemma}}
\def\bdf{\begin{defi}}
\def\edf{\end{defi}}
\def\btm{\begin{theorem}}
\def\etm{\end{theorem}}
\def\bex{\begin{exe}}
\def\eex{\end{exe}}
\def\bpp{\begin{propos}}
\def\epp{\end{propos}}
\def\bQ{\begin{qes}}
\def\eQ{\end{qes}}
\def\ben{\begin{enumerate}}
\def\een{\end{enumerate}}
\def\into{\hookrightarrow}

\def\t{\tilde}
\def\lb{\label}

\newcommand{\qbin}[2]{\bigl[{\begin{matrix}\raise 0.32ex \hbox{\mathsmall{10}{#1}}\vspace*{-2.5mm}  \\ \mathsmall{10}{#2}\end{matrix}}\bigr]}
\newcommand{\qBin}[2]{\Bigl[{\begin{matrix}\mathsmall{12}{#1}\vspace*{-1.2mm} \\ \mathsmall{12}{#2}\end{matrix}}\Bigr]}

\newcommand{\Q}{{\mathbb Q}}

\begin{document}

\institute{}

\title{
The Lawrence-Krammer-Bigelow Representations of the Braid Groups  via $U_q({\mathfrak s}{\mathfrak l}_2)$}

\author{Craig Jackson  \ \ and \ \  Thomas Kerler}
 
\maketitle

%%%%%%%%%%%   Begin -- Chapter 1 -- Introduction %%%%%%%%%%%%%%%%%%%%%%

\begin{center}

\today

\begin{abstract} We construct representations of the braid groups $B_n$ on $n$ strands on  free ${\mathbb Z}[q^{\pm 1},s^{\pm 1}]$-modules
$\W_{n,l}$ using generic Verma modules for an integral version of $\uqsl$. We prove that the $\W_{n,2}$ are  isomorphic to the faithful
 Lawrence Krammer Bigelow  representations of $B_n$ after appropriate identification of parameters of Laurent polynomial rings by constructing
explicit integral bases and isomorphism. We also prove that the $B_n$-representations 
$\W_{n,l}$ are irreducible over the fractional
field $\Q (q,s)$.  

\end{abstract}

\end{center}

%%%%%%%%%%%%%%%

\section{Introduction}
In recent years 
the representation theory of the braid groups $\,B_n\,$ on $n$ strands 
has attracted attention due to two groundbreaking developments. 
One of them is in the work of Bigelow and 
Krammer \cite{Big01, Kra02}, who managed to resolve the long standing
problem of the linearity of the braid groups by showing that a two-parameter 
generalization of the classical Burau representation is faithful. The second
development is
the emergence of  vast families of braid group representations 
that are constructed  from  quantum algebras (see \cite{KRT97} and references therein)
and conformal field theories \cite{KZ84,Koh87,TK88}. 
Intriguing relationships 
between these seemingly very different approaches 
have been discovered and they remain a fascinating area of study.

\smallskip

In this article we give an  explicit construction and proof of 
an isomorphism between the faithful representation $\lkb{2}{n}$ of $B_n$  
considered by Bigelow and Krammer and the submodule of  
the {\em R-matrix} representations on $\Vm^{\otimes n}$ for the generic Verma
module  $\Vm$ of the quantum group $\,U_q({\mathfrak s}{\mathfrak l}_2)\,$. 

%##############################################3

For the purpose of this article we will consider Krammer's version 
 $\lkb{2}{n}$, as defined in \cite{Kra02} and restated in \cite{Big03}, 
and refer to it as the 
 {\em Lawrence-Krammer-Bigelow representation} or LKB representation. 
It  is defined over the 
ring $\,{\mathbb Z}[{\tt q}^{\pm 1},{\tt t}^{\pm 1}]\,$ of two variable 
Laurent polynomials with integral coefficients. The parameters
${\tt q}$ and ${\tt t}$ are associated to Deck transformations of 
a covering $\tilde C_n\to C_n$, where $C_n$ is the two-point configuration
space on a disc with  $n$-punctures. The natural representation of $B_n$
on $H_2(\tilde C_n)$ as a  $\,{\mathbb Z}[{\tt q}^{\pm 1},{\tt t}^{\pm 1}]\,$  
is isomorphic to $\lkb{2}{n}$ over $\mathbb Q({\tt q}^{\pm 1},{\tt t}^{\pm 1})\,$, see  \cite{Big01}.
While these modules are not isomorphic over 
$\,{\mathbb Z}[{\tt q}^{\pm 1},{\tt t}^{\pm 1}]\,$  (see \cite{PP02}), Bigelow conjectures in \cite{Big03} that
the relative homology $H_2(\tilde C_n,\tilde \nu)$ is isomorphic to 
$\lkb{2}{n}$
over $\,{\mathbb Z}[{\tt q}^{\pm 1},{\tt t}^{\pm 1}]$ , where $\tilde \nu$ may be understood as
a piece in the boundary of a certain compactification of $\tilde C_n\,$.

\smallskip

The first obstacle in finding such an isomorphism is that 
the braid group representations obtained from quantum groups are originally
defined over the complex numbers rather than integral two-variable Laurent 
polynomials.  To this end we will define $\,U_q({\mathfrak s}{\mathfrak l}_2)\,$ 
as an algebra $\Alg$ over ${\mathbb Z}[q,q^{-1}]$, and introduce the generic 
Verma module $\Vm$ over  $\,\Rqs={\mathbb Z}[q^{\pm 1},s^{\pm 1}]\,$, 
where $s$ may be 
thought of as the exponential highest weight  $\,s=q^{\lambda}\,$. 
The braid group action  commutes with the $\Alg$-action so that the 
highest weight spaces
$\, \W_{n,l} \subset {\bf V}^{\otimes n} \,$   of  $\Alg$, corresponding to 
weights $sq^{-2l}=q^{\lambda-2l}$, are again $B_n$-invariant. 
We prove in Section~\ref{sec-IHW} that the $\,\W_{n,l}\,$ are free  $\,\Rqs$-modules, and 
construct explicit bases
\begin{equation}\label{eq-basis}
\Bas_{n,l}=\{w_{\vec\alpha} \mid \vec\alpha=(\alpha_j,\ldots,\alpha_{n})\ 
     \mbox{ with } j > 1 \mbox{ and } \sum_{i=j}^n\alpha_i=l-1\}\,
\end{equation}
such that $\,\W_{n,l}\,$ is the $\,\Rqs$-span of $\;\Bas_{n,l}$. Specifically, we find

\medskip

\begin{theorem}\label{thm-freemodule}
The highest weight space $\,\W_{n,l}\subset{\bf V}^{\otimes n}\,$ is a free module over 
$\,\Rqs={\mathbb Z}[q^{\pm 1},s^{\pm 1}]\,$ with explicitly given 
basis $\,\Bas_{n,l}\,$ as in (\ref{eq-basis}).  Hence, for each $\,l\geq 0\,$  we
obtain a  representation of the braid group $B_n$ in $\,n$-strands given 
by a homomorphism as follows:
\beq
\rho_{n,l}\,:\;B_n\,\longrightarrow \,{\rm GL}\!
\left({\mbox{\small
$\binom{n+l-2}{l}$
},\Rqs
}\right)\;\cong\;{\rm Aut}_{\RqsFN}(\W_{n,l})
\eeq
\end{theorem}

The identification of the quantum representation on $\,\W_{n,2}\,$ 
from Theorem~\ref{thm-freemodule} with the  LKB representation $\lkb{2}{n}$ further 
requires an identification of parameters which we give by the 
following  monomorphism between Laurent polynomials. 
\begin{equation}\label{eq-deftheta}
\theta \,:\;{\mathbb Z}[{\tt q}^{\pm 1},{\tt t}^{\pm 1}]\,\longrightarrow\,{\mathbb Z}[{s}^{\pm 1},{q}^{\pm 1}]=\Rqs\;:\;\;
\begin{array}{rcl}
{\tt q} & \mapsto & s^{2}\\
{\tt t}&\mapsto & -q^{-2}\rule{0mm}{6mm}
\end{array}
\end{equation}
Consider also the involutive automorphism $\invo$ of $B_n$ defined on the
generators by $\invo(\sigma_i)=\sigma^{-1}_i\,$ (given by switching all crossing or reflection at 
the plane of projection of a braid), and denote by $\,\lkb k n ^{\dagger}$ the 
representation given by pre-composing the action on $\,\lkb k n$ with $\invo$.
With these conventions the main result of this article, which we will prove in Section~\ref{sec-lkb}, can be formulated as follows:

\medskip

\begin{theorem}\label{thm-main}
For every $n>1$ there is a  isomorphism of $B_n$-representations over $\Rqs$ 
\beq\label{eq-mainiso}
\W_{n,2}\quad\stackrel{\cong}{\longrightarrow}\quad\lkb{2}{n}^{\dagger}\otimes_{\theta}\Rqs
\eeq
which maps the basis $\Bas_{n,2}$ to the fork basis from \cite{Big01}.
\end{theorem}

In \cite{Zin00} Zinno manages to find a different  identification of  the 
LKB representation with  a quantum algebraic object,
namely the quotient  of the Birman-Wenzl-Murakami algebra  similarly defined 
over $\,{\mathbb Z}[{\tt q}^{\pm 1},{\tt t}^{\pm 1}]\,$.  
This representation can, by \cite{Wen90}, be understood as the one 
arising from the quantum orthogonal groups 
$\,U_{\zeta}({\mathfrak s}{\mathfrak o}(k+1))\,$ acting on 
the $\,n$-fold tensor product of the fundamental
representation. Since the representation in  \cite{Wen90} 
is irreducible this 
implies that $\lkb{2}{n}$ and hence $\W_{n,2}$ are
irreducible for all $n>1$.  
\footnote{Note, however, that in the symmetric group specialization 
with  $s=1$ and $q=1$ these
representation are clearly reducible for all $l\geq 2$}

In Section~\ref{sec-Irred} we generalize this result in our case to obtain

\btm\label{thm-irred}
For all $n \geq 2$ and $l \geq 0$ the $B_n$-representation $\W_{n,l}$
is irreducible over the fraction field $\overline{\Rqs} = \Q (q,s)$.
\etm

\medskip

Faithfulness of $\W_{n,l}$ for $l\geq 3$ is still an open question, 
as are identifications of these
representation with geometrically constructed ones analogous to 
Theorem~\ref{thm-main}. Obvious
candidates for a generalization of  Theorem~\ref{thm-main} 
are the $B_n$-representations constructed by Lawrence in 
\cite{Law90}. The starting point there is again the  configuration 
space $\,Y_{n,l}\,$ of $\,l\,$ points in the plane with $\,n\,$ holes.
The braid group action is then naturally defined on 
$\lkb l n =H_l(\tilde Y_{n,l})$, where $\tilde Y_{n,l}$ is the canonical
cover of $Y_{n,l}$ with covering group ${\mathbb Z}^2$.  
The latter makes the representation spaces into 
$\,{\mathbb Z}[{\tt q}^{\pm 1},{\tt t}^{\pm 1}]$-modules.
\medskip

\begin{conjecture}\label{conj-lawrence-iso} The spaces  
$H_l(\tilde Y_{n,l})$ are free 
$\,{\mathbb Z}[{\tt q}^{\pm 1},{\tt t}^{\pm 1}]$-modules
which carry  an (irreducible) action of  $B_n$ as defined in \cite{Law90}.
They are isomorphic to the representations of $B_n$ on weight spaces $\W_{n,l}$ 
over $\Rqs$ after appropriate identifications of parameters in
the Laurent polynomial rings. 
\end{conjecture}

The first obvious piece evidence for this conjecture is that it holds 
for $l\leq 2$. Indeed, for $l=1$ both $\lkb 1 n =H_1(\tilde Y_{n,1})$ and 
$\W_{n,1}$ can be readily identified with the classical Burau 
representation of $B_n$.  For more details see the beginning of Section~\ref{sec-lkb}.

It has been observed, both by Lawrence (Section~4 of  \cite{Law90}) and by Bigelow
(Section~6 of \cite{Big03}), that for $l=2$ and the parameter specialization ${\tt t}=-{\tt q}^{-1}$
the  LKB-representation has as a factor the Temperley-Lieb representation associated to the two-row
partition $[n-2,2]$. In the former case it occurs as a quotient and in the latter setting
as a sub-module.

In Section~\ref{sec-TLS} we will explain the occurrence of the Temperley-Lieb factor from the 
point of view of quantum-$\sl$ representations. Particularly, the respective identification
$s=q$ will correspond to specializing the highest weight of the fundamental representation
of quantum-$\sl$ within the Verma module. The exact sequence of $B_n$-modules we establish
in (\ref{eq-sesBn}) reflects the cohomological picture of \cite{Law90}.

In Lemma~\ref{lm-TLinLKB}  we will also identify the irreducible $n$-dimensional quotient of the specialized 
LKB-representation by the Temperley-Lieb sub-representation, and prove that the sequence in  (\ref{eq-sesBn}) 
does not split. Consequently, although the representation becomes reducible in the ${\tt t}=-{\tt q}^{-1}$
specialization it remains indecomposable. We also discuss in Section~\ref{sec-TLS} the construction of 
braid elements in the kernel of the Temperley-Lieb representation in order to underscore the loss of 
information in the parameter specialization.

Theorem~\ref{thm-main} as well as its generalization in   
Conjecture~\ref{conj-lawrence-iso} are inspired by \cite{Fel91} and \cite{VS91} where 
quantum-$\,{\mathfrak s}{\mathfrak l}_2$ actions on the homology of local 
systems over similar configuration spaces are constructed.

\medskip

\noindent
{\bf Acknowledgments:} The second author thanks Giovanni Felder for 
very useful discussions about  \cite{Fel91} which motivated this article. 
We also thank the anonymous referees for numerous suggestions that 
helped to improve  the article and led to the addition of Section~\ref{sec-TLS}.

%%%%%%%%%%%   End -- Chapter 1 -- Introduction %%%%%%%%%%%%%%%%%%%%%%
%%%%%%%%%%%%%  Begin -- Chapter 2   %%%%%%%%%%%%%%%%%%%%%%%%%%%%%%%%%

\section{From Topological to Integral Braid Group Representation} \label{S-tensor}

In this section we review the basic definitions and constructions of quantum  $\sl$ which lead to the relevant 
representations of the braid groups. We will start from the framework of quasi-triangular topological Hopf algebras
due to Drinfeld  \cite{Dri87} over rings of power series. An exposition and further development of Drinfeld's theory
can be found in Kassel's textbook \cite{Kas95} which we will use as main reference.

We start with the definition of the algebra $\Alg_{\hbar}$ over a power series ring $\mathbb P[[\hbar]]$ 
where $\mathbb P$ is some commutative ring containing the rational numbers $\mathbb Q$. 
The indeterminate is related to $h$ used in  \cite{Kas95} by $\hbar=\frac 12 h\,$.
The generators of $\Alg_{\hbar}$ are $E$, $F$, and $H$  with relations
\begin{equation}\label{eq-Hrels}
\begin{array}{rcr}
\,[H,E] &=& 2E\\
\,[H,F] &=& -2F\rule{0mm}{7mm}
\end{array}
\qquad\quad
[E,F] = \frac {\sinh(  \hbar H)} {\sinh(  \hbar )}\;\;.
\end{equation}
The algebra $\Alg_{\hbar}$ is given by formal power series $\sum_na_n\hbar^n$ where each coefficient
$a_n$ is a finite combination of monomials in the generators $E$, $F$, and $H$ over $\mathbb P\supseteq \mathbb Q$.
It is easy to see that the expression for $[E,F]$ can indeed be written in this way.
In addition, the comultiplication on $\Alg_{\hbar}$ is defined by
\begin{equation}\label{eq-Hcorels}
\begin{array}{rcl}
\Delta(E)& = &E\otimes e^{  \hbar H}\,+\,1\otimes E\\
 \Delta(F)& = &F\otimes 1\,+\,e^{-  \hbar H}\otimes F\rule{0mm}{7mm}
\end{array},\quad\mbox{and}\quad
\Delta(H)=H\otimes 1\,+\,1\otimes H
\end{equation}
Formally, the coproduct is a homomorphism $\Delta: \Alg_{\hbar}\to\Alg_{\hbar}\widetilde\otimes \Alg_{\hbar}$, where the 
tensor completion is described in  Section~XVI.3 of \cite{Kas95}. We introduce the usual set of notations for 
$q$-numbers, $q$-factorials, and $q$-binomial coefficients:
\begin{equation}
\begin{array}{rcl}
q\,&=&e^{  {\hbar} }\\
 \,[n]_q&=&\displaystyle \frac {q^n-q^{-n}}{q-q^{-1}}=\frac {\sinh\left({ {\hbar}  n}\right)} {\sinh\left({  {\hbar}  }\right)}\rule{0mm}{8mm}
\end{array}
\qquad
\begin{array}{rcl}
\,[n]_q!&=&[n]_q[n-1]_q\ldots[2]_q[1]_q\\
\qbin{n}{j}_q&=&\displaystyle \frac {[n]_q!}{[n-j]_q![j]_q!}\rule{0mm}{8mm}
\end{array}
\end{equation} 

Note that all of these quantities are invertible in $\mathbb P[[\hbar]]$ for $n\neq 0\,$.
A {\em universal R-matrix} for $\Alg_{\hbar}$ is now given as in Theorem~XVII.4.2 of \cite{Kas95}  by
\begin{equation}
\label{eq-univRmat}
{\mathcal R}\;=\; 
e^{\frac {\hbar}2(H\otimes H)} \cdot \Bigl (
\sum_{n=0}^{\infty}q^{\frac {n(n-1)}2}\frac {(q-q^{-1})^n}{[n]_q!}
E^n\otimes F^n \Bigr )\quad\in\; \Alg_{\hbar}\widetilde\otimes \Alg_{\hbar}
\end{equation} 

Drinfeld's construction  from   \cite{Dri87} as described in \cite{Kas95} implies that the {\em R}-matrix from 
(\ref{eq-univRmat}) makes $\Alg_{\hbar}$ into a quasi-triangular topological Hopf algebra. Particularly, this implies 
that  ${\mathcal R}$ obeys the {\em Yang-Baxter} relation given as an equation in 
$\Alg_{\hbar}\widetilde\otimes \Alg_{\hbar}\widetilde\otimes \Alg_{\hbar}$ by 
\begin{equation}\label{eq-YBR}
 {\mathcal R}_{12}{\mathcal R}_{13}{\mathcal R}_{23}\;=\;{\mathcal R}_{23}{\mathcal R}_{13}{\mathcal R}_{12}\;.
\end{equation}
Moreover, $ {\mathcal R}$ fulfills the usual commutation relation in $\Alg_{\hbar}\widetilde\otimes \Alg_{\hbar}$ given by 
\begin{equation}\label{eq-YBcomm}
 {\mathcal R}\Delta(x)\,=\,\Delta^{opp}(x){\mathcal R}\qquad \forall x\in \Alg_{\hbar}\;.
\end{equation}

In order to construct representations of the braid groups we will need to consider first representations
of $\Alg_{\hbar}$. Instead of distinguishing many representations by their highest weights we consider only one representation and ``absorb'' the highest weight as a parameter in the underlying 
coefficient ring as follows.

In  \cite{Kas95} the coefficient ring was chosen as $\mathbb P=\mathbb C\,$, yet all calculations and
statements there clearly also apply for any other choice of $\mathbb P\supseteq \mathbb Q\,$.  For our purposes 
we will choose the coefficient ring to be $\mathbb P =\mathbb Q[\lambda]$, that is,  the polynomial ring with 
rational coefficients in one indeterminate $\lambda$ which may be thought of as a generic highest weight.

$\Alg_{\hbar}$ is thus an algebra over  $\mathbb Q[\lambda][[\hbar]]$ -- the ring of power series in $\hbar$ 
whose  coefficients are rational polynomials in $\lambda$. In this setting $\Alg_{\hbar}$ admits a special highest 
weight module over the same ring described as follows.

Consider  the $\mathbb Q[\lambda]$-module $\Ql$ freely generated by an infinite sequence of vectors denoted by  $\,\{v_0,v_1,\ldots\}\,$.
The {\em generic Verma module} $\Vmh=\Ql[[\hbar]]$ is then the associate topologically free module by  in the sense of 
Section~XVI.2 of \cite{Kas95}. The action of $\Alg_{\hbar}$ on $\Vmh$ is given by 
\begin{eqnarray}\label{eq-Haction}
 H.v_j\,&=&\,(\lambda-2j)v_j\nonumber \\
E.v_j\,&=&\,v_{j-1}\\
F.v_j\,&=&\,[j+1]_q\cdot[\lambda-j]_q\,v_{j+1}\nonumber
\end{eqnarray}

Note here that indeed $\,[\lambda-j]_q=\frac {\sinh(\hbar(\lambda-j))} {\sinh(\hbar)}
\in\,\mathbb Q[\lambda][[\hbar]]\,$. The module is similar to the standard highest 
weight module obtained as the
induced representation associated to the one-dimensional representation of the Borel algebra generated
by $E$ and $H$ acting on $v_0$. It is, however, not equivalent to this module since the elements
$[\lambda-j]_q$ are not invertible in $\,\mathbb Q[\lambda][[\hbar]]\,$. (Ring evaluations $\lambda =m$ create
additional highest weight vectors for the traditional Verma module, but additional lowest weight
vectors for the representation in (\ref{eq-Haction}). The modules are equivalent, and irreducible, only
for evaluations $\lambda\not\in\mathbb N$.) 

As described in the end of Sections~XVI.4 any topological $\Alg_{\hbar}$-module such as $\Vmh$ now entails a 
solution to the Yang-Baxter equation on $\Vmh\widetilde\otimes \Vmh\widetilde\otimes \Vmh$ by (\ref{eq-YBR}), 
which  commutes with the action of $\Alg_{\hbar}$ on the same space by (\ref{eq-YBcomm}). As an endomorphism 
on $\Vmh\widetilde\otimes \Vmh$ we define the action of a braid group generator by 
\begin{equation}\label{eq-BraidRepT}
\Rm\,:\; \Vmh^{\widetilde\otimes 2}\,\rightarrow\,\Vmh^{\widetilde\otimes 2}\,:\qquad  v\widetilde \otimes w\;\mapsto\;
e^{-\frac{\hbar}2\lambda^2}T(\mathcal R.(v\widetilde \otimes w))
\end{equation}

Here $\mathcal R$ acts as an element of  $\Alg_{\hbar}\widetilde\otimes\Alg_{\hbar}$ on $\Vmh\widetilde\otimes \Vmh$, 
and $T$ denotes the usual transposition $T(v\otimes w)=w\otimes v$. We also multiply the map by the unit
$e^{-\frac{\hbar}2\lambda^2}\,\in\,\mathbb Q[\lambda][[\hbar]]\,$ which also yields a solution to the Yang Baxter equation  
since this relation is homogeneous. 

The braid group $B_n=\langle \sigma_1\ldots\sigma_{n-1}\,|\,\sigma_i\sigma_{i+1}\sigma_i=\sigma_{i+1}\sigma_i
\sigma_{i+1},\; \sigma_i\sigma_j=\sigma_j\sigma_i \mbox{ for }|i-j|>1\rangle$ is now represented on 
$\Vmh^{\widetilde\otimes n}$  by the assignment
\begin{equation}\label{eq-sigmatilde}
 \sigma_i\quad\mapsto\quad \widetilde\sigma_i=\Bid^{\widetilde\otimes i-1}\widetilde\otimes 
\,\Rm\, \widetilde\otimes 
\Bid^{\widetilde\otimes n-i-1}
\end{equation}

The goal of the following constructions is to identify a  sublattice in $\Vmh^{\widetilde\otimes n}\,$
which is  invariant under this action of the braid group $B_n$. This lattice will be a free module over
a subring $\,\Rqs\subset \mathbb Q[\lambda][[\hbar]]\,$ which is characterized as follows. 
Consider first the following ring homomorphism from the two-variable Laurent polynomials to the 
power series ring: 
\begin{equation}\label{eq-laurent2power}
i_{\hbar}\,:\;\Rqs =\mathbb Z[q,q^{-1},s,s^{-1}]\,\longrightarrow\,\mathbb Q[\lambda][[\hbar]]\;:\quad 
\begin{array}{rcl}
q&\mapsto & e^{\hbar}\\
s&\mapsto & e^{\hbar\lambda}\;.
\end{array}
\end{equation}
It is clear that $i_{\hbar}$ is well defined by inspection of the power series expansion in $\hbar\,$, and that 
$i_{\hbar}$ is a monomorphism since $(\hbar,\lambda)\mapsto (e^{\hbar},e^{\hbar\lambda})$ has dense
image in $\mathbb C^2$. We will thus denote the image of  $i_{\hbar}$  
also by $\,\Rqs=\mathbb Z[q^{\pm 1}, s^{\pm 1}]\subset \,\mathbb Q[\lambda][[\hbar]]\,$ with the identification
of parameters as prescribed in (\ref{eq-laurent2power}).

In order to find a suitable subalgebra over this ring we define next a set of special generators in $\Alg_{\hbar}$ by
\begin{equation}
K\;=\;e^{\hbar H}\,,\qquad K^{-1}\;=\;e^{-\hbar H}\,,\qquad\mbox{and} \qquad F^{(n)}\;=\;\mathsmall{10}{\displaystyle \frac {(q-q^{-1})^n}{[n]_q!}}F^n\;.
\end{equation}
The generators $ F^{(n)}$ are similar to the divided powers introduced by Lusztig in \cite{Lus93} but differ by the 
additional $(q-q^{-1})$ factors.
The following relations readily follow from the ones given in (\ref{eq-Hrels}):
\begin{equation}\label{eq-Urels}
\begin{split}
KK^{-1}=K^{-1}K=1\,,\quad\quad 
KEK^{-1}&=q^{2}E\,,\quad\quad 
KF^{(n)}K^{-1}\;=\;q^{-2n}F^{(n)}
\\
F^{(n)}F^{(m)}\;=\;\qBin{n+m}{n}_q F^{(n+m)}\,,
\quad&{\rm and}\quad
[E, F^{(n+1)}]\;=\;F^{(n)} (q^{-n}K-q^{n}K^{-1})\;.
\end{split}
\end{equation}

Let now $\Alg\subset\Alg_{\hbar}$ be the subalgebra over $\Rqs$ generated by the set of
elements $\,\{K,E,F^{(n)}\}$. As a sublattice $\Alg$ is the free $\Rqs$-module spanned
by the PBW basis $\,\{K^lE^mF^{(n)}\,:\,l\in\mathbb Z, m,n\in\mathbb N\cup\{0\}\}\,$. In fact,   
 $\Alg$ is isomorphic to the algebra over $\Rqs$ defined abstractly by  generators 
$\,\{K^{\pm 1}, E, F^{(n)}\}\, $ and the relations given in (\ref{eq-Urels}).

The coproduct and antipode evaluated on the generators of $\Alg$ are readily computed:
\begin{equation} \label{eq-Ucop}
\begin{split}
\Delta(K)\,=\,K\otimes K\qquad & \qquad \Delta(E)\,=\,E\otimes K+1\otimes E\\
\Delta(F^{(n)})\;=&\;\sum_{j=0}^n q^{-j(n-j)}K^{j-n}F^{(j)}\otimes F^{(n-j)}\\
S(K)\;=\;K^{-1}\hspace*{12mm}S(E)\;&=\;-EK^{-1}\hspace*{12mm}S(F^{(n)})\;=\;(-1)^nq^{n(n-1)}K^{n}F^{(n)}
\end{split}
\end{equation}
These formulas immediately imply that the coproduct is in fact a map $\,\Delta:\,\Alg\to\Alg\otimes\Alg\,$
with $\otimes$ taken over $\Rqs$. Consequently, $\Alg$ is a Hopf subalgebra 
of $\Alg_{\hbar}$, and thus a Hopf algebra over $\Rqs$ by itself. 

\medskip

Next, let  $\,\Vm\subset \Vmh\,$ be the free $\Rqs$-module generated
by the basis vectors $\,\{v_0,v_1,\ldots\}\,$. That is, an element in $\Vm$ is given by 
$\sum_j p_jv_j$ with $p_j\in \mathbb \Rqs=Z[q^{\pm 1}, s^{\pm 1}]$ and only finitely many $p_j$ 
are non-zero. 
The actions of the generators of  $\Alg$ on the basis vectors $v_j$ are easily worked out from the action of 
$\Alg_{\hbar}$ to be the following: 
\begin{equation}\label{NewVermaEq1}
\begin{split}
K . v_j = sq^{-2j} v_j \hspace*{3mm} & \hspace*{12mm}
E . v_j =  v_{j-1}\\
F^{(n)} . v_j &= \left({
   \qBin{n+j}{j}_q
   \prod_{k=0}^{n-1} (sq^{-k-j}-s^{-1}q^{k+j})}\right) v_{j+n}.
\end{split}
\end{equation} 

Observe that all coefficients in these formulae lie in the subring  $\Rqs=\mathbb Z[q^{\pm 1}, s^{\pm 1}]$
and contain only a finite number (one) of vectors. This immediately implies the following:

\begin{lemma}\label{eq-U-Lat}
The subspace $\Vm\subset \Vmh$ is invariant under the action of the subalgebra $\Alg\subset \Alg_{\hbar}$.
\end{lemma}

This also implies that the natural
actions of $\Alg$ as well as $\Alg^{\otimes n}$ on $\Vmh^{\widetilde\otimes n}$ map the respective
subspace $\Vm^{\otimes n}\subset \Vmh^{\widetilde\otimes n}\,$ to itself. The main observation of this
section is that the same is true for the braid group action.

\begin{lemma}\label{lm-R-Lat}
The map $\Rm\,$, as defined in (\ref{eq-BraidRepT}), maps the subspace $\,\Vm^{\otimes 2}\subset\Vmh^{\widetilde\otimes 2}$ to itself.
\end{lemma}

\begin{proof} We first note that the map $\Rm$ can be written as the composite of three maps
\begin{equation}\label{eq-R=TCP}
\Rm\;=\;T\circ {\mathcal C}\circ{\mathcal P}\;,
\end{equation}
where $T$ is the usual transposition as in (\ref{eq-BraidRepT}). The operator $\mathcal C$ is given by the action of the
factor $e^{\frac {\hbar}2(H\otimes H)}\,$ from the expression in (\ref{eq-univRmat}) for the universal $R$-matrix 
multiplied by the extra term $e^{-\frac {\hbar} 2\lambda^2}$ that occurs in (\ref{eq-BraidRepT}). Finally, $\mathcal P$
is given by application of the remaining summation in parentheses in    (\ref{eq-univRmat}).

We prove that each of these three operators in (\ref{eq-R=TCP}) preserves $\,\Vm^{\otimes 2}\,$ as a subspace. 
This is trivially true for $T$.  For the action of  $\mathcal C$ we compute
\begin{eqnarray}
{\mathcal C}.(v_j\otimes v_k)\,&=&\,e^{-\frac {\hbar} 2\lambda^2}e^{\frac {\hbar}2(H\otimes H)}v_j\otimes v_k
\,=\,e^{-\frac {\hbar} 2\lambda^2}e^{\frac {\hbar}2(\lambda-2j)(\lambda-2k)}v_j\otimes v_k \nonumber \\
\,&=&\,e^{ -\hbar \lambda(j+k)+2\hbar j k  }v_j\otimes v_k
\,=\,s^{-(j+k)} q^{2jk}v_j\otimes v_k\rule{0mm}{7mm}
\end{eqnarray}
Thus ${\mathcal C}.(v_j\otimes v_k)\in \Vm^{\otimes 2}\,$ and the claim follows for $\mathcal C$. 

For $\mathcal P$ we first rewrite the summation expression for the universal $R$-matrix in  (\ref{eq-univRmat})
in terms of the generators of $\Alg$.
\begin{equation}\label{eq-Pexp}
{\mathcal P}\;=\;\sum_{n=0}^{\infty}q^{\frac {n(n-1)}2}\frac {(q-q^{-1})^n}{[n]_q!}
E^n\otimes F^n\,
\;=\;
\sum_{n=0}^{\infty}q^{\frac {n(n-1)}2}
E^{n}\otimes F^{(n)}\;.
\end{equation}
The fact that the action of $E$ on $\Vm$ is locally nilpotent together with the observation that any finite 
truncation of the summation in (\ref{eq-Pexp}) yields an element  in $\Alg\otimes\Alg$  imply the
claim for $\mathcal P$. More specifically, the action of $\mathcal P$ can be worked out explicitly
to be the following. 
\begin{equation}
{\mathcal P}.(v_i\otimes v_j)\;=\;\sum_{n=0}^iq^{\frac {n(n-1)}2} 
   \qBin{n+j}{j}_q
   \prod_{k=0}^{n-1} (sq^{-k-j}-s^{-1}q^{k+j}) \,v_{i-n}\otimes v_{j+n}.
\end{equation}
Since the summation is a finite one and all coefficients are in  $\mathbb Z[q^{\pm 1}, s^{\pm 1}]$ we
can now infer that ${\mathcal P}.(v_i\otimes v_j)\in\Vm^{\otimes 2}\,$. Consequently, all three operators
$T$, $\mathcal C$, and $\mathcal P$ map the subspace $\Vm^{\otimes 2}\,$ to itself, which proves the
lemma.
\end{proof}

\medskip

For future use let us also record the explicit formula for the action of $\Rm$ on $\,\Vm^{\otimes 2}\,$. 
\begin{equation}\label{eq-intRmat}
\begin{split}
 \Rm.(v_i\otimes v_j)&\;=\; \\
s^{-(i+j)}\sum_{n=0}^i  & q^{2(i-n)(j+n)}q^{\frac {n(n-1)}2} 
   \qBin{n+j}{j}_q
   \prod_{k=0}^{n-1} (sq^{-k-j}-s^{-1}q^{k+j}) \,v_{j+n} \otimes v_{i-n}\,.
\end{split}
\end{equation}

Let us summarize our finding of this section in the following theorem:

\begin{theorem}\label{thm-braid-rep}
The maps $\, \sigma_i= \Bid^{\otimes i-1}\otimes \,\Rm\,\otimes \Bid^{\otimes n-i-1}\,$, 
with the 
$\,\Rm\,$ as in (\ref{eq-intRmat}), define a representation of the braid group $B_n$ 
on $\,\Vm^{\otimes n}\,$, as a free $\,\mathbb Z[s^{\pm 1},q^{\pm 1}]$-module. 
The maps $\sigma_i$ also commute with
the action of $\Alg$ on $\,\Vm^{\otimes n}\,$ and preserve the natural $\mathbb Z$ grading.
\end{theorem}

\begin{proof}
The fact that the maps $\sigma_i$ preserve $\,\Vm^{\otimes n}\,$ is immediate from Lemma~\ref{lm-R-Lat}.
They fulfill the braid group relations since they are restrictions of the maps $\widetilde\sigma_i$
from (\ref{eq-sigmatilde}) which fulfill these relations by construction. Moreover, these maps commute
with the action of $\Alg_{\hbar}$ and hence also with the action of $\Alg$. 
\end{proof}

%%%%%%%%%%%%%  End -- Chapter 2   %%%%%%%%%%%%%%%%%%%
%%%%%%%%%%%%%  Begin -- Chapter 3   %%%%%%%%%%%%%%%%%

\section{Integrality of Highest Weight Spaces}\label{sec-IHW}

The main purpose of this section is to prove the assertion 
in Theorem~\ref{thm-freemodule}, namely, that the highest weight spaces are
free $\Rqs$-modules.

In order to define these highest weight spaces 
let  $\Vm_{n,l} = \ker (K - s^n q^{-2l}) \subset {\Vm}^{\otimes n}$ 
be the weight space corresponding to the weight $s^n q^{-2l}$.
Recall that $x \in \Alg$ acts on ${\Vm}^{\otimes n}$ 
by $\Delta^{(n)} x$, where $\Delta^{(n)}: \Alg \to \Alg^{\otimes n}$ is defined recursively by $\Delta^{(2)} = \Delta$ and  $\Delta^{(n)} = (\Delta^{(n-1)} \otimes 1 ) \Delta$.
By \eqref{eq-Ucop} and \eqref{NewVermaEq1}, $\Vm_{n,l}$ is the $\Rqs$-span of the vectors
$v_{\alpha_1} \otimes \cdots \otimes  v_{\alpha_n}$ where 
$\alpha_1 + \cdots + \alpha_n = l$. We now define 
\begin{equation}\label{eq-defW}
 \W_{n,l} = \ker (E) \cap \Vm_{n,l}.
\end{equation}
The space
$\W_{n,l}$ is the so-called {\em highest weight space}
corresponding to the weight $s^n q^{-2l}$.
Since the representation of
$B_n$ on $\Vm^{\otimes n}$ commutes with the $\Alg$-action,
we see that both  $\Vm_{n,l}$ and $\W_{n,l}$ are also 
$B_n$ representations. 

Let us also define $\Am_{n,l}, \Bm_{n,l} \subset \Vm_{n,l}$ for $l \geq 2$ by
\begin{equation}\label{eq-Adef}
\begin{split}
 \Am_{n,l} &= \text{$\Rqs$-span of }  \ABas_{n,l}\hfill   \\
\text{with} & \quad \ABas_{n,l} =\{
 v_{\alpha_1} \otimes \cdots \otimes  v_{\alpha_n} 
  \mid \exists k \text{ such that } \alpha_k = 1 
       \text{ and } \alpha_j = 0 \text{ } \forall \text{ }  j<k \}
\end{split}
\end{equation}
and
\begin{equation}\label{eq-Bdef}
\begin{split}
\Bm_{n,l} &= \text{$\Rqs$-span of }  \BBas_{n,l}\hfill   \\
\text{with} & \quad \BBas_{n,l} =\{ v_{\alpha_1} \otimes \cdots \otimes v_{\alpha_n}
  \mid \exists k \text{ such that } \alpha_k > 1
       \text{ and } \alpha_j = 0 \text{ } \forall \text{ }  j<k \}.
\end{split}
\end{equation}

We immediately see that $\Vm_{n,l} = \Am_{n,l} \oplus \Bm_{n,l}$. 
Given a multi-index $\vec \alpha = (\alpha_{j}, \dots, \alpha_n )$
for some $j > 1$ such that
$\sum_{i=j}^n \alpha_i = l-1$, we can define an element of $\ABas_{n,l}$  by
\beq\label{Abasis}
a_{\vec{\alpha}} = v_0^{\otimes (j-2)} \otimes v_1 \otimes 
    v_{\vec{\alpha}}
\eeq 
where $v_{\vec{\alpha}} = 
    v_{\alpha_{j}} \otimes \cdots \otimes v_{\alpha_n}
    \in \Vm_{n-j+1, \, l-1}$.
Clearly, letting $\vec \alpha$ vary among all such multi-indices
gives the basis $\ABas_{n,l}$ of $\Am_{n,l}$.

\blm\label{EonB}
For all $n \geq 1$ and $l \geq 2$, the map 
$E|_{\sBm_{n,l}} : \Bm_{n,l} \to \Vm_{n,l-1}$ is an $\Rqs$-linear isomorphism.
\elm
\begin{proof}
To show that $E|_{\sBm_{n,l}} : \Bm_{n,l} \to \Vm_{n,l-1}$ is surjective
we need to show that for every $v_{\vec{\alpha}} = v_{\alpha_1} 
\otimes \cdots \otimes v_{\alpha_n} \in  \Vm_{n,l-1}$ there is 
some $b \in \Bm_{n,l}$ such that $E.b = v_{\vec{\alpha}}$. 

We proceed by induction
on $j = l- \alpha_k$, where $\alpha_k$ is the first nonzero entry in the
multi-index $\vec{\alpha} = (\alpha_1, \dots, \alpha_n )$. The initial case, 
when $j=1$, occurs when $\alpha_k = l-1$ and is handled simply 
by observing that $E . v_l = v_{l-1}$.
To prove the induction step let us take $v_{\vec{\alpha}} =
v_0^{\otimes (k-1)} \otimes v_{\alpha_k}
\otimes \cdots \otimes v_{\alpha_n}
\in 
V_{n,l-1}$ such that
$l - \alpha_k = j+1$. Setting
$b = v_0^{\otimes (k-1)} \otimes v_{\alpha_k + 1}
\otimes \cdots \otimes v_{\alpha_n}$
we see that $b \in \Bm_{n,l}$ and 
\beq
E.b = (\text{unit})\, v_{\vec{\alpha}} + (\text{other terms})
\eeq
where the first nonzero index in each of the 
other terms is $\alpha_k +1$. Hence, the other 
terms satisfy the induction hypothesis and so are in the image of $E$. From this
it follows that $v_{\vec{\alpha}}$ is in the image of $E$.

To show that $E|_{\sBm_{n,l}}$ has no kernel take some $0 \neq b \in \Bm_{n,l}$.
Then $b$ will have some minimal
term in its expression, namely, some
$v_{\alpha_1} \otimes \cdots \otimes v_{\alpha_k} \otimes
  \cdots \otimes v_{\alpha_n}$ where $\alpha_i = 0$ for all $i < k$,
$\alpha_k \geq 2$, and if $v_{\beta_1} \otimes \cdots 
 \otimes v_{\beta_n}$
is in the expression for $b$ then $\beta_i = 0$ for all $i<k$ and either
$\beta_k =0$ or $\beta_k \geq \alpha_k$. Then comparing the terms in
the expression for $E.b$ we see that it is impossible to
cancel out the term
$v_{\alpha_1} \otimes \cdots \otimes v_{\alpha_k -1} \otimes
  \cdots \otimes v_{\alpha_n}$. Hence, $E.b \neq 0$.
\end{proof}

Since $E |_{B_{n,l}}$ is an isomorphism, we seek a way to parametrize $\mbox{Ker}(E)$ by $A_{n,l}$. This parametrization is accomplished with an $\Rqs$-linear map $\Phi : \Vm_{n,l} \to \Vm_{n,l}$, constructed in such a way that $E \circ \Phi$ vanishes on $A_{n,l}$ (see Lemma \ref{lm-EPhi-form}). Hence, for $l \geq 2$, define $\Phi$ on basis elements 
$ a_{\vec{\alpha}}=v_0^{\otimes (j-2)} \otimes v_1 \otimes 
              v_{\vec{\alpha}}\,\in \ABas_{n,l}$ and $ b\in \BBas_{n,l}$ as follows:
\begin{equation}\label{map}
\begin{split}
 \Phi (a_{\vec{\alpha}})  &  = 
  \sum_{k = 0}^l b_{\vec \alpha, k} \,
     v_0^{\otimes (j-2)} \otimes v_{k} \otimes 
             E^{k-1} v_{\vec{\alpha}}        \\
\Phi (b)  &=  b  \;.  \rule{0mm}{7mm}
\end{split}
\end{equation}
The coefficients are given by 
\beq\label{c}
b_{\vec \alpha , k} = (-1)^{k-1} s^{(k-1)(j-n-1)} q^{(k-1)(2l - k -2)}. 
\eeq
Notice that when $k=0$ in \eqref{map} we have a multiple of the term 
$v_0^{\otimes (j-1)} \otimes E^{-1} v_{\vec{\alpha}}$. 
By $E^{-1} v_{\vec{\alpha}}$ we mean the unique element 
$\eta \in \Bm_{n-j+1,l}$ such that $E \eta =  v_{\vec{\alpha}}$.
Such an element $\eta$ exists and is unique because of Lemma~\ref{EonB}.

\begin{lemma}\label{lm-PhiIso}
We have $(\Phi-\Bid)^2=0$ so that $\Phi$ is an automorphism of $\Vm_{n,l}$.
\end{lemma}

\begin{proof}
Clearly, we have $(\Phi-\Bid)(b)=0$ for $b\in\Bm_{n,l}$.
 For $k=1$ we have $b_{\vec \alpha, 1} = 1$ so that
$(\Phi-\Bid)(a_{\vec \alpha})=\Phi (a_{\vec \alpha}) - a_{\vec \alpha} \in \Bm_{n,l}$ and hence
$(\Phi-\Bid)^2(a_{\vec \alpha})=0$. The nilpotency relation immediately implies that
$\Phi^{-1}=2-\Phi$ is an inverse.
\end{proof}

Under the change of basis on $\Vm_{n,l}$ given by $\Phi$ the operator $E$ has a simple form.

\begin{lemma}\label{lm-EPhi-form} For all $n \geq 1$ and $l \geq 2$ the composite $E\circ\Phi$ vanishes on $\Am_{n,l}$ and is injective on $\Bm_{n,l}$ with
\begin{equation}\label{eq-decEPhi}
 E\circ \Phi\,=\,0\,\oplus\,E|_{\sBm_{n,l}}\,:\;\Am_{n,l}\,\oplus\,\Bm_{n,l}\,\longrightarrow \,\Vm_{n,l-1}\;.
\end{equation}
This implies that the following is an isomorphism of $\Rqs$-modules:
\begin{equation}\label{eq-PhiIsoAW}
\Phi\,:\,\Am_{n,l}\,\stackrel{\cong}{\longrightarrow}\,\W_{n,l}\;.
\end{equation}
\end{lemma}

\begin{proof} The first half of the action in (\ref{eq-decEPhi}) is to show $E\circ\Phi$ is zero on any element
$a_{\vec{\alpha}}$ which we verify by explicit computation:
\begin{align*}
E \circ \Phi(a_{\vec{\alpha}}) 
 &= {\Delta}^{(n)}(E). \sum_{k \geq 0} b_{\vec \alpha, k} \,
     v_0^{\otimes (j-2)}\! \otimes v_{k} \otimes
             E^{k-1} v_{\vec{\alpha}}\\
 &= \sum_{k \geq 1} s^{n-j+1} q^{-2(l-k)} b_{\vec \alpha, k}  \,
     v_0^{\otimes (j-2)}\! \otimes v_{k-1} \otimes
             E^{k-1} v_{\vec{\alpha}}  +
    \sum_{k \geq 0} b_{\vec \alpha, k} \,
     v_0^{\otimes (j-2)}\! \otimes v_{k} \otimes
             E^{k} v_{\vec{\alpha}}\\
 &= \sum_{k \geq 0} \big( s^{n-j+1} q^{-2(l-k-1)} b_{\vec \alpha, k+1} 
                      + b_{\vec \alpha, k} \big) \,
     v_0^{\otimes (j-2)} \otimes v_{k} \otimes
             E^{k} v_{\vec{\alpha}} \;  =  \;0 \;.
\end{align*}
Here we use that \eqref{c} implies the recursion
\[
s^{n-j+1}q^{-2(l-k-1)} b_{\vec \alpha, k+1}
= - b_{\vec \alpha, k}\;.
\]
Now (\ref{eq-decEPhi}) follows from   (\ref{map}), where $\Phi$ is defined to be identity on $\Bm_{n,l}\,$.
By Lemma~\ref{EonB} we have that $E|_{\sBm_{n,l}}$ is injective so that 
$\,\mathrm{ker}(E\circ \Phi)\cap \Vm_{n,l}\,=\,\Am_{n,l}\,$. Since, by Lemma~\ref{lm-PhiIso}, $\Phi$ is
an automorphism of $\Rqs$-modules this implies (\ref{eq-PhiIsoAW}).
\end{proof}

Let us also describe the case $l=1$ more explicitly. A basis of 
$\Vm_{n,1}$ is given by 
\begin{equation}\label{eq-c-bas}
 c_i\,=\,v_0^{\otimes (i-1)}\otimes v_1\otimes v_0^{\otimes (n-i)}\qquad \mbox{for }
i=1,\ldots, n
\end{equation}
The subspaces defined in (\ref{eq-Adef}) and (\ref{eq-Bdef}) are defined slightly
different for $l=1$, namely
\begin{equation}\label{eq-ABl=1}
\begin{split}
\Am_{n,1}\,& =\,\Rqs\mbox{-span of }\ABas_{n,1}\qquad \mbox{with}\quad \ABas_{n,1}=\{c_i|i=1,\ldots,n-1\}\\
\Bm_{n,1}\,& =\,\Rqs\mbox{-span of }\BBas_{n,1}\;\qquad \mbox{with}\quad \BBas_{n,1}=\{c_n\}\;.\\
\end{split}
\end{equation}
In this setting we have 
$E^{-1}(v_0^{\otimes m})=v_0^{\otimes (m-1)}\otimes v_1\,\in\Bm_{m,1}$. Formula
(\ref{map}) thus yields a basis for $\W_{n,1}$ given by vectors
\begin{equation}\label{eq-defwl=1}
w_i\,=\,\Phi(c_i)\,=\,c_i-s^{(n-i)}c_n\;\qquad \mbox{with } i=1,\ldots, n-1\;.
\end{equation}
With these conventions it is easy to see that all previous lemmas in this section
also apply to the case $l=1$ (and trivially so to the case $l=0$). 

\smallskip 

We are now in a position to prove Theorem \ref{thm-freemodule}, namely, that the 
highest weight spaces are free $\Rqs$-modules.

\begin{proof}[Proof of Theorem \ref{thm-freemodule}] Since (\ref{eq-PhiIsoAW}) is an isomorphism of
$\Rqs$-modules and $\Am_{n,l}$ is clearly a free module, also $\W_{n,l}$ has to be a free $\Rqs$-module. 
The rank is given by the number of vectors in the set of spanning vectors given in (\ref{eq-Adef}), which 
is given by $\binom{n+l-2}{l}$.
\end{proof}

Since the generators $\sigma_i\,$, as defined in (\ref{eq-sigmatilde}), map (by $\Alg$-equivariance) 
each $\W_{n,l}$ subspace to itself, Lemma~\ref{lm-EPhi-form} implies that the conjugate maps 
$\sigma_j^{\Phi}=\Phi^{-1}\circ \sigma_i\circ \Phi$ map $\Am_{n,l}$ to itself. Thus the representation of
$B_n$ over $\Rqs$ given by $\sigma_j|_{\sW_{n,l}}$ is equivalent to the representation given by the maps
$\sigma_j^{\Phi}|_{\sAm_{n,l}}\,$.

Suppose  $\pi_{A}$ is the 
projection of $\Vm_{n,l}$ onto $\Am_{n,l}$ along $\Bm_{n,l}$.  Observe also that 
$\,\Phi^{-1}|_{\W_{n,l}}\,=\,\pi_A|_{\W_{n,l}}$. This yields the basic but useful formula:
\begin{equation}\label{basiscalc}
 \sigma_j^{\Phi}|_{\sAm_{n,l}}\;=\;\pi_A\circ\sigma_j\circ\Phi\;.
\end{equation}
Implicit to this formula is the method of calculating the action of a braid generator $\sigma_j^{\Phi}$ on a 
particular basis vector:
\begin{enumerate}
\item[{\sf 1)}] For a basis vector $\,a_{\vec \alpha}\in \ABas_{n,l}\,$ determine $\Phi(a_{\vec \alpha})\in\W_{n,l}$ by (\ref{map}).
\item[{\sf 2)}] Use (\ref{eq-intRmat}) and (\ref{eq-sigmatilde}) to determine the image $\sigma_j(\Phi(a_{\vec \alpha}))$.
\item[{\sf 3)}] Write $\sigma_j(\Phi(a_{\vec \alpha}))$ in the standard basis $\ABas_{n,l}\cup \BBas_{n,l}$ 
and eliminate the components of $\BBas_{n,l}\,$ leaving an $\Rqs$-linear combination of vectors from $\ABas_{n,l}$.
\end{enumerate}

In the following we also consider the action of $B_n$ directly on $\W_{n,l}\subset \Vm_{n,l}$. 
A natural basis is given by  $\Bas_{n,l}\,=\,\Phi(\ABas_{n,l})\,=\,\{w_{\vec{\alpha}} = \Phi(a_{\vec{\alpha}})\}\,$. 
By construction the explicit action of the braid 
group generators $\sigma_j$ in this basis is exactly the same as the action
of the $\sigma_j^{\Phi}$ in the basis $\ABas_{n,l}$ so that the computations
remain the same. 

%%%%%%%%%%%%%  End -- Chapter 3   %%%%%%%%%%%%%%%%%%
%%%%%%%%%%%%%  Begin-- Chapter 4   %%%%%%%%%%%%%%%%%

\section{The Lawrence Krammer Bigelow Representation}\label{sec-lkb}

Here we prove that the representation of $B_n$ on $\W_{n,2}$ is isomorphic the LKB representation
which was recently shown in \cite{Big01} and \cite{Kra02} to be faithful. 
As preparation let us show first
 that the representation $\W_{n,1}$ is isomorphic to the classical, reduced Burau
representation over $\mathbb Z[{\tt t},{\tt t}^{-1}]\,$.

The formula for the $R$-matrix in (\ref{eq-intRmat}) implies 
$\Rm.(v_0\otimes v_0)=v_0\otimes v_0\,$, $\Rm.(v_0\otimes v_1)=s^{-1}v_1\otimes v_0\,$, and $\Rm.(v_1\otimes v_0)=s^{-1}v_0\otimes v_1+(1-s^{-2})v_1\otimes v_0\,$.
Applied to the basis $\{c_i\}$ of $\Vm_{n,1}$ from (\ref{eq-c-bas}) this implies
the following action of $B_n$ on $\,\Vm_{n,1}\,$:
\begin{equation}\label{eq-sbur}
\begin{split}
\sigma_i . c_j &= c_j \phantom{space} j \neq i, i+1\\
\sigma_i . c_i &= s^{-1} c_{i+1} + (1-s^{-2}) c_i\\
\sigma_i . c_{i+1} &= s^{-1} c_i.
\end{split}
\end{equation}
Using the rescaled basis $\{d_j=s^{-j}c_j\,|\,1\leq i <n\}$ and with a substitution of parameter $\,s^{-2}\mapsto {\tt t}$ the action from (\ref{eq-sbur}) turns out 
to yield
exactly the {\em unreduced} Burau representation  $\lkbt 1 n$ of dimension $n$ as described, 
for example, in (3-23) of \cite{Bir74}. Thus we have by identification of basis vectors that 
\begin{equation}
\Vm_{n,1} \cong \lkbt 1 n\otimes_{{\tt t}=s^{-2}}\Rqs\;.
\end{equation}
Now, the basis for $\W_{n,1}$ from (\ref{eq-defwl=1}) may also be rescaled as
\begin{equation}\label{eq-ubas}
u_j\,=\,s^{j}w_j\,=\,s^{2j}d_j-s^{2n}d_n\,=\,{\tt t}^{-j}d_j-{\tt t}^{-n}d_n \qquad \mbox{with } j=1,\ldots, n-1\;.
\end{equation}

Recall that the {\em reduced} Burau representation $\lkb 1 n$ of dimension $(n-1)$  is given by the kernel of the map $\,\lkbt 1 n \to \mathbb Z[{\tt t}^{\pm 1}]:\,d_j\mapsto {\tt t}^j\,$. Clearly, the basis described in (\ref{eq-ubas}) is thus a basis also of $\lkb 1 n$ and we obtain the following relation.
\begin{lemma} \qquad $\displaystyle \W_{n,1} \cong \lkb 1 n\otimes_{{\tt t}=s^{-2}}\Rqs\;.$
\end{lemma}

Let us now turn to the $l=2$ case. The basis  $\ABas_{n,2}$ from \eqref{Abasis} is given by elements
\begin{equation}\label{eq-defaij}
 a_{i,j} =   v_0^{\otimes (i-1)} \otimes   v_1 \otimes
            v_0^{\otimes (j-i-1)} \otimes   v_1 \otimes
           v_0^{\otimes (n-j)} \qquad\quad\mbox{for } 1 \leq i < j \leq n\;.
\end{equation}
Correspondingly, the $\BBas_{n,2}$ consists of the following elements:
\begin{equation}\label{eq-defbk}
 b_k =   v_0^{\otimes (k-1)} \otimes   v_2 
                      \otimes   v_0^{\otimes (n-k)}\qquad\quad\mbox{for } 1 \leq k \leq n\;.
\end{equation}

The basis $\Bas_{n,2}$
for $\W_{n,2}$ is  given by application of the map in \eqref{map}  to
$\ABas_{n,2}$ which yields the following set of elements:
\begin{equation}\label{eq-defwij}
w_{i,j} = \Phi (a_{i,j}) = a_{i,j} - s^{j-i} q^{-2} b_j - s^{i-j} b_i \qquad\quad\mbox{for } 1 \leq i < j \leq n\;.
\end{equation}
 
The action of the braid group $B_n$ on these vectors is now 
computed using the step by step procedure following 
(\ref{basiscalc}).  In addition to the formulae in the previous paragraph this also involves calculating
expressions for $\Rm.(v_i\otimes v_j)\,$ with $i+j=2$. 

In each of these expressions only the coefficients of the $v_1\otimes v_1$-term needs to be considered
since the contributions of the $v_0\otimes v_2$-terms and $v_2\otimes v_0$-terms will be projected out by $\pi_A\,$. 
The relevant relations are thus the following: 
\begin{equation}\label{eq-R2ord}
\begin{array}{lll}
\Rm.(v_0\otimes v_2)\,&=\, 0 &\,\mod \langle v_2\otimes v_0,v_0\otimes v_2\rangle \;\\
\Rm.(v_1\otimes v_1)\,&=\,q^2 s^{-2}(v_1\otimes v_1) 
&\,\mod \langle v_2\otimes v_0,v_0\otimes v_2\rangle \; \\
\Rm.(v_2\otimes v_0)\,&=\,q^2 (s^{-1}-s^{-3})(v_1\otimes v_1) 
&\,\mod \langle v_2\otimes v_0,v_0\otimes v_2\rangle \; \\
\end{array}
\end{equation}
Applying (\ref{eq-R2ord}) to the elements in (\ref{eq-defaij}) and (\ref{eq-defbk}), and 
combining expressions in (\ref{eq-defwij}) we can compute the action of $B_n$ on the
basis vectors in $\Bas_{n,2}$
according to the procedure given at the end of the previous section. 
The resulting formulae for the generators of $B_n$ are listed next where we assume that
$\{i, i+1 \} \cap \{j,k \} = \emptyset$:
\begin{equation}\label{eq-sigma.w}
\begin{array}{ll}
\sigma_i . w_{j, k}  &=  w_{j,k} \arh \\
\sigma_i . w_{i+1, j}   &= s^{-1} w_{i,j} \arh \\
\sigma_i . w_{j, i+1}  &= s^{-1} w_{j,i}  \arh\\
\sigma_i . w_{i, j}  &=  s^{-1} w_{i+1, j} + (1-s^{-2}) w_{i,j} - s^{i-j-1} (1-s^{-2}) q^2 w_{i,i+1} \arh \\
\sigma_i . w_{i, i+1}  &= s^{-4} q^2 w_{i, i+1} \arh \\
\sigma_i . w_{j, i}   &=     s^{-1} w_{j,i+1} + (1-s^{-2}) w_{j, i} - s^{i-j-1}(1-s^{-2}) w_{i,i+1}.  \arh
\end{array}
\end{equation}

For comparison we consider the explicit Lawrence Krammer Bigelow representation 
$\lkb 2 n $ of $B_n$ as given in Section~5.2 of
\cite{Big03}. (Note that the representation given in \cite{Big01} contains a sign error which is corrected in \cite{Big03}). There the space $\lkb 2 n $ is described as the free $\mathbb Z[{\tt t}^{\pm 1},{\tt q}^{\pm 1}]$-module spanned by basis elements $\,\{F_{i,j}:1\leq i < j \leq n\}\,$. From the formulae in 
\cite{Big03} the actions of the inverses of the generators of $B_n$ on $\lkb 2 n $ 
are readily worked out to be as follows:
\begin{equation}\label{eq-sigma.F}
 \begin{array}{ll}
\sigma_i^{-1}. F_{j, k}     &=  F_{j,k}  \arh \\
\sigma_i^{-1}. F_{i+1, j}  &=  F_{i,j} \arh \\
\sigma_i^{-1}. F_{j, i+1}  &=  F_{j,i}  \arh \\
\sigma_i^{-1}. F_{i, j}      &=  {\tt q}^{-1} F_{i+1, j} + (1-{\tt q}^{-1}) F_{i,j} +{\tt t}^{-1} ({\tt q}^{-1}-{\tt q}^{-2})  F_{i,i+1} \arh \\
\sigma_i^{-1}. F_{i, i+1}  &= -{\tt t}^{-1} {\tt q}^{-2} F_{i, i+1}  \arh\\
\sigma_i^{-1}. F_{j, i}     &= {\tt q}^{-1} F_{j,i+1}  + (1-{\tt q}^{-1}) F_{j, i}  - ({\tt q}^{-1}-{\tt q}^{-2}) F_{i,i+1}. \arh \\
\end{array}
\end{equation}

\begin{proof}[Proof of Theorem~\ref{thm-main}]
Let us define the map $\mathscr{F}\,:\,\lkb 2 n \to \W_{n,2}\,$ by
\[
\mathscr{F}(F_{i,j})=s^{i+j}w_{i,j}\qquad \mbox{and}\qquad \mathscr{F}(pv+qw)=\theta(p)\mathscr{F}(v)+
\theta(q)\mathscr{F}(w) \;,
\] 
where $p,q\in \mathbb Z[{\tt t}^{\pm 1},{\tt q}^{\pm 1}]$, $\;v,w\in\lkb 2 n $, and  
$\theta$ is the ring homomorphism given in 
(\ref{eq-deftheta}). It follows now by direct computation from the equations in
(\ref{eq-sigma.w}) and (\ref{eq-sigma.F}) that
\[
 \mathscr{F}\sigma_i^{-1}\;=\;\sigma_i\mathscr{F} \quad \forall 
i \in \{1, \dots, n \}
\qquad\mbox{so that}\qquad 
\mathscr{F}\invo(\beta) \;=\; \beta \mathscr{F} \qquad {\forall \beta \in B_n}
\]
where $\invo$ is the involution described in the introduction. Hence $\,\mathscr{F}:\,\lkb 2 n^{\dagger} \to \W_{n,2}\,$ is $B_n$-equivariant by definition.
Since it also maps basis vectors to basis vectors of free modules and $\theta$ is
a monomorphism, $\lkb 2 n ^{\dagger}$ can be considered a $B_n$-submodule of $\W_{n,2}$
whose  $\Rqs$-span is again $\W_{n,2}$. This implies the isomorphism in
Theorem~\ref{thm-main}.
\end{proof}

%%%%%%%%%%%%%  End -- Chapter 4   %%%%%%%%%%%%%%%%%%%%
%%%%%%%%%%%%%  Start  -- Chapter 5   %%%%%%%%%%%%%%%%%

\section{The Temperley-Lieb Specialization}\label{sec-TLS}

 In Section~6 of \cite{Big03} Bigelow considers the parameter specialization  ${\tt q}{\tt t}=-1$
for a version of the LKB-representation, and recovers  a sub-module on which the $B_n$-action
factors through the respective Temperley-Lieb algebra with representation associated to a
two-row Young tableau. The latter, in turn, are closely related to  the representation theory of
quantum-$\sl$ via Schur-Weyl duality.

In this section we will show how the Temperley-Lieb  submodule structure  naturally follows by extracting 
finite-dimensional highest or lowest weight modules of quantum-$\sl$ from the generic Verma modules 
used in Theorem~\ref{thm-braid-rep} for respective parameter identifications in the ground ring.. The Temperley-Lieb algebra
then arises as the centralizer in the case of the tensor powers of the 2-dimensional fundamental representation
of quantum-$\sl$.

The topological and representation theoretic derivations of the same submodule structure in Theorem~6.1 of 
\cite{Big03} and  Lemma~\ref{lm-TLinLKB}  below, respectively, give thus another insight into the topological 
content of quantum-$\sl$ actions. In addition, we will address in Lemma~\ref{lm-TLinLKB} the splitting property
and complementary module structure, and conclude with general remarks on the loss of information in the
Temperley-Lieb reduction.

In order to construct finite dimensional quantum-$\sl$ representations we
 fix a positive integer $\ell\in\mathbb N$ and consider the module with ring quotient
into $\mathbb Z[s^{\pm 1},q^{\pm 1}]\to\mathbb Z[q,q^{-1}]$ that sends $s\mapsto q^{\ell}$. This yields 
$\Alg$-modules over $\mathbb Z[q,q^{-1}]$ defined as follows:
\begin{equation}\label{eq-Vmm}
\breve\Vm_{\ell}=\Vm \!\!\underset{s \,=\, q^{\ell}}{\otimes}\!\mathbb Z[q,q^{-1}]\;.
\end{equation}

Clearly, $\breve\Vm_{\ell}$ is still a free $\mathbb Z[q,q^{-1}]$-module with basis $\{v_0, v_1, \ldots\}$. 
It is immediate from (\ref{NewVermaEq1}) that
\begin{equation}\label{eq-FonVmm}
 F^{(n)}.v_j\,=\,0 \qquad \mbox{ for } \; j+n > \ell
\end{equation}
Suppose $\Im_{\ell}\subset\breve\Vm_{\ell}$ is the free $\mathbb Z[q,q^{-1}]$-submodule spanned by $\{v_0, v_1, \ldots, v_{\ell}\}\,$. 
It follows easily from (\ref{eq-Vmm}) and  (\ref{NewVermaEq1}) that $\Im_{\ell}$ is also a $\Alg$-submodule, that is,
$\Alg.\Im_{\ell}=\Im_{\ell}$. It may be thought of as the irreducible lowest weight module whose lowest weight vector $v_{\ell}$
has the properties $\,K.v_{\ell}=q^{-\ell}v_{\ell}\,$ and $\,F^{(n)}.v_{\ell}=0\,$ for $n\geq 1\,$.
It also follows readily, for example from (\ref{eq-intRmat}), that 
\begin{equation}\label{eq-RinvIl}
\,\Rm.(\Im_{\ell}\otimes \Im_{\ell})\subseteq \Im_{\ell}\otimes \Im_{\ell}\;.
\end{equation}
Thus we can specialize and restrict the braid group representations from Theorem~\ref{thm-braid-rep} to the following
finite rank module over $\mathbb Z[q,q^{-1}]$.
 \begin{equation}
 \Im_{\ell}^{\otimes n}\quad \subseteq \quad  \breve\Vm_{\ell}^{\otimes n} \quad =\quad \Vm^{\otimes n}\!\!\underset{s \,=\, q^{\ell}}{\otimes}\!\mathbb Z[q,q^{-1}]\;.
 \end{equation}
These braid group representations are equivalent over $\mathbb Q(q)$ to the ones obtained from 
the standard $R$-matrix construction for the $(\ell+1)$-dimensional
representations of quantum-$\mathfrak{sl}_2$ (for example, Section~VIII.3 in \cite{Kas95}), 
and also correspond to limits of solutions to the Yang-Baxter equation given in \cite{KRS81}. 
As before, the highest weight constructions yield respective sub-representations of the braid groups.
 For the following discussion let us instead consider all relevant modules over the 
field of fractions $\mathbb Q(q)$:
\begin{equation}\label{eq-WLmm}
\begin{split}
 \breve\Vm_{n,k,\ell}=\Vm_{n,k} \!\!\underset{s \,=\, q^{\ell}}{\otimes}\!\mathbb Q(q) &\qquad\mbox{and}\quad 
\breve\W_{n,k,\ell}=\W_{n,k} \!\!\underset{s \,=\, q^{\ell}}{\otimes}\!\mathbb Q(q)\\
\mathbf{L}_{n,k,\ell}=\breve\W_{n,k,\ell}\cap \Imo_{\ell}^{\otimes n}   &
 \qquad\;\mbox{with}\quad  \Imo_{\ell}=\Im_{\ell}{\otimes}\mathbb Q(q)
\end{split}
% \qquad\;\mbox{with}\quad 
% \quad\mbox{and}\quad 
\end{equation}

Of particular interest is the specialization $\ell=1$, that is, $s=q$, which corresponds to the fundamental representation of quantum-$\mathfrak{sl}_2$.
In this case 
$\,\Imo_1\,=\,\mathbb Q(q)v_0\,\oplus \,\mathbb Q(q)v_1\,$ so that the $R$-matrix acts on a
4-dimensional space spanned by $v_0\otimes v_0$, $v_1\otimes v_0$, $v_0\otimes v_1$, 
and $v_1\otimes v_1\,$. The action is more conveniently described in terms of
\begin{equation}
 \mathsf{E}\,:=\,q(\Rm-\Bid^{\otimes 2})
\end{equation}
for which we can compute readily from the explicit formula (\ref{eq-intRmat}) that
\begin{equation}\label{eq-Eact}
\begin{split}
 \mathsf{E}(v_0\otimes v_0)\,&=\,0\,=\,\mathsf{E}(v_1\otimes v_1)\\
\mathsf{E}(v_0\otimes v_1) \,&=\, v_1\otimes v_0 \,-\, q(v_0\otimes v_1)  \\
\mathsf{E}(v_1\otimes v_0)\,&=\, v_0\otimes v_1 \,- \, q^{-1} (v_1\otimes v_0)\\
\end{split}
\end{equation}
The formulae in (\ref{eq-Eact}) can, in turn, be used to verify the following relations:
\begin{equation}
\begin{split}
 \mathsf{E}^2 \,&=\,-(q+q^{-1}) \mathsf{E}\\
 (\mathsf{E}\otimes \Bid) (\Bid\otimes\mathsf{E} ) (\mathsf{E}\otimes \Bid)\,& =\,\mathsf{E}\otimes \Bid\\
  (\Bid\otimes\mathsf{E} )(\mathsf{E}\otimes \Bid) (\Bid\otimes\mathsf{E} )\,& =\,\Bid\otimes\mathsf{E}\\
\end{split}
\end{equation}
These relations are easily recognized as those of the Temperley-Lieb algebra  $\mathsf{A}_{n,q}\,$. 

Over the fraction field $\mathbb Q(q)$ (or over $\mathbb C$ with $q$ specialized to a value that is not a 
root of unity) it is well known that the images of  $\mathsf{A}_{n,q}\,$ and $\Alg$ in $\mathrm{End}(\Imo_1^{\otimes n})$
via the obvious representations are semisimple and each others commutants, see \cite{Ji86}. This implies the quantum
analogue of Schur-Weyl duality, namely that the $n$-fold tensor product is isomorphic over $\mathbb Q(q)$ to
\begin{equation}\label{eq-SWdec}
 \Imo_1^{\otimes n} \;\cong\;\sum_{k=0}^{\lfloor \frac n 2 \rfloor}F_{[n-k,k]}\otimes  {\pi}_{[n-k,k]}
\end{equation}
as a $\Alg \times \mathsf{A}_{n,q}\,$-module. Here $F_{[n-k,k]}$ is the representation of highest
weight $q^{(n-2k)}$, and ${\pi}_{[n-k,k]}$ the $\,\mathsf{A}_{n,q}\,$ representation associated to the partition
$[n-k,k]$ in analogy to the symmetric group \cite{Jon87}. 
The dimensions of the factors are the same as in the classical theory (see for example Section~9 of \cite{GW09}):
\begin{equation}\label{eq-dimTL}
 \dim( {\pi}_{[n-k,k]})={n\choose k} -{n\choose k-1}  \qquad\mbox{ and} \qquad \dim(F_{[n-k,k]}) =n+1-2k.
\end{equation}
Suppose  $v_0^{k}$ is the highest weight vector of  $F_{[n-k,k]}$. It follows readily from (\ref{eq-SWdec}) that the space of
highest weight vectors of weight $q^{(n-2k)}$corresponds to $\langle v_0^{k}\rangle\otimes  {\pi}_{[n-k,k]}$. Thus 
with definitions from (\ref{eq-defW}) and (\ref{eq-WLmm}) we obtain the following identification of  $\,\mathsf{A}_{n,q}$-modules:
\begin{equation}\label{eq-piL}
 {\pi}_{[n-k,k]}\;\cong\;\mathbf{L}_{n,k,1}\;\subseteq\;\breve\W_{n,k,1}\;.
\end{equation}
In order to apply this to the situation of the LKB representation let us denote by $\tau$  the following ring
homomorphism
\begin{equation}\label{eq-deftau}
\tau \,:\;{\mathbb Z}[{\tt q}^{\pm 1},{\tt t}^{\pm 1}]\,\longrightarrow\,{\mathbb Q}(q)\;:\;\;
\begin{array}{rcl}
{\tt q} & \mapsto & q^{2}\\
{\tt t}&\mapsto & -q^{-2}\rule{0mm}{6mm}
\end{array}\;
\end{equation}
We also introduce an $n$-dimensional representation $\mathbf C_n(\lambda)$. To this end, let $B_n\to \mathbb Z$ be the Abelian quotient map (with $\sigma_i \mapsto1$) 
and $B_n\to S_n:\,b\,\mapsto\overline b$ the symmetric group quotient. Then let  
$\mathbf C_n(\lambda)=\langle e_1,\ldots, e_n\rangle$ where elements of $B_n$ acts as
\begin{equation}\label{def-eqCact}
 \sigma_j.e_j= \lambda e_{j+1}\,,\quad \sigma_j.e_{j+1}= e_{j}, \quad\mbox{and}\quad  \sigma_j.e_{i}=e_i\;\;\mbox{for } i\not\in\{j,j+1\}\;.
\end{equation}
We can now state the following  relation of the LKB representation with the Temperley-Lieb representation theory.
\begin{lemma}\label{lm-TLinLKB}
Reducing the ground ring of the  LKB representation by $\tau$ to $\mathbb Q(q)$ as in (\ref{eq-deftau})
we obtain for $n\geq 4$ the 
following short exact sequence of $B_n$-modules
\begin{equation}\label{eq-sesBn}
0\,\to\, {\pi}_{[n-2,2]}\;\hookrightarrow\;
 \lkb{2}{n}^{\dagger}\!\underset{\tau}{\otimes}\mathbb Q(q)\;\twoheadrightarrow\;\mathbf C_n(q^{-4})\;\to\;0\;
\end{equation}
where the $\mathbb Q(q)[B_n]$-action on the first summand factors through $\,\mathsf{A}_{n,q}\,$and the $B_n$-action 
on the 
second through the combined quotient $\mathbb Z\times S_n$.
For $n\geq 4$ the sequence in (\ref{eq-sesBn}) is {\em not} split.
\end{lemma}

\begin{proof} The  inclusion given the second map in (\ref{eq-sesBn}) is the same as the inclusion in 
(\ref{eq-piL}) via the identifications  $\mathbf{\pi}_{[n-2,2]}\overset{(\ref{eq-piL})} {\cong}\mathbf{L}_{n,2,1}
\overset{(\ref{eq-WLmm})}{=}\breve\W_{n,2,1}\cap\Imo_1^{\otimes n}\hookrightarrow \breve\W_{n,2,1} 
\overset{(\ref{eq-WLmm})}{=}\W_{n,2}\underset{q=s}{\otimes}\mathbb Q(q)
\overset{(\ref{eq-mainiso})}{\cong}(\lkb{2}{n}^{\dagger}\!\underset{\theta}{\otimes}\mathbb Z[s^{\pm 1},q^{\pm 1}])\underset{q=s}{\otimes}\mathbb Q(q)
=\lkb{2}{n}^{\dagger}\!\underset{\tau}{\otimes}\mathbb Q(q)\,
$. The cokernel of this inclusion naturally maps to the following quotient  of weight spaces:
\begin{equation}\label{eq-WVquot}
 J\,:\;
\frac {\breve\W_{n,2,1}}{\breve\W_{n,2,1}\cap \Imo_{1}^{\otimes n}}\;\longrightarrow\;
\frac {\breve\Vm_{n,2,1}}{\breve\Vm_{n,2,1}\cap \Imo_{1}^{\otimes n}}
\end{equation}
A basis over $\mathbb Q(q)$ of $\breve\Vm_{n,2,1}$ is given by the $\ABas_{n,2}=\{a_{i,j}\}_{1\leq i < j\leq n}$ and
$\BBas_{n,2}=\{b_{k}\}_{1\leq k\leq n}$as defined in (\ref{eq-defaij}) and (\ref{eq-defbk}). Clearly, the subspace 
$\breve\Vm_{n,2,1}\cap \Imo_{1}^{\otimes n}$ is exactly the subspace spanned by $\ABas_{n,2}$ so that the
quotient on the right side of (\ref{eq-WVquot}) is an $n$-dimensional space for which a basis is given by the classes
$\overline b_k$ of the basis elements $b_k$.

It follows readily that $J(\overline w_{i,j})=-q^{j-i-2}\overline b_j- q^{i-j}\overline b_i$ where the
generators $\overline w_{i,j}$ are the respective classes of the basis elements from  (\ref{eq-defwij}).
It is a straightforward exercise in linear algebra  to show that every $\overline b_k$ can be expressed
as a combination of elements $-q^{j-i-2}\overline b_j- q^{i-j}\overline b_i$ over $\mathbb Q(q)$ if $n\geq 3\,$.

Thus $J$ is a surjective map. Using the fact that the LKB representation has dimension $n\choose 2$, formulae (\ref{eq-dimTL}) and \eqref{eq-piL} together
imply that the domain of $J$ is also $n$-dimensional. Consequently,
$J$ is an isomorphism of $B_n$-modules. 

The module structure on the image of $J$ is found by computing the action of the $R$-matrix on 
$\breve\Vm_{n,2,1}$ modulo $\Imo_1^{\otimes n}$. Specializing $q=s$ in
(\ref{eq-intRmat}) we  find
\begin{equation}
\begin{split}
\Rm(v_2\otimes v_0)\;& =\;q^{-2}v_0\otimes v_2\;+\;(q-q^{-1})v_1\otimes v_1\\
\Rm(v_0\otimes v_2)\;& =\;q^{-2}v_2\otimes v_0 \quad\mbox{and}\quad 
\Rm(v_0\otimes v_0)\;=\;v_0\otimes v_0\;.
\end{split}
\end{equation}
Thus, the respective action on  $\overline b_k\in\breve\Vm_{n,2}/\Imo_1^{\otimes n}$ is given by 
\begin{equation}
 \sigma_k.\overline b_k=q^{-2} \overline b_{k+1} \qquad
\sigma_k.\overline b_{k+1}=q^{-2} \overline b_{k}\qquad 
\sigma_k.\overline b_j=\overline b_{j}  \mbox{ for } j\not\in\{k,k+1\} \;.
\end{equation}
Upon setting $\lambda=q^{-4}$ and 
after renormalization of the basis 
\begin{equation}
 e_j \,= \, - q^{2j}\overline b_j 
\end{equation}
this is precisely the same action as the one described in (\ref{def-eqCact}), and hence proves the exact 
sequence in (\ref{eq-sesBn}).

In order to show that this sequence is not split for $n\geq 4$ it suffices to show that 
$0\to \breve\Vm_{n,2,1}\cap \Imo_{1}^{\otimes n}\to
\breve\Vm_{n,2,1}\to \mathbf C_n(q^{-4})\to 0$ is not split since any splitting homomorphism for (\ref{eq-sesBn})
can be composed with the inclusion $\breve\W_{n,2,1}\hookrightarrow\breve\Vm_{n,2,1}$. Such a splitting would
imply the existence of generators $e_j\in \breve\Vm_{n,2,1}$ for $j=1,\ldots, n$ with a $B_n$ action as prescribed
in   (\ref{def-eqCact}) for $\lambda=q^{-4}$ and with $\,e_j\equiv -q^{2j}\overline b_j\,\mod \Imo_1^{\otimes n}$.

The minimal polynomial of $\sigma_i$ on $\breve\Vm_{n,2,1}$ is given by 
$\mu(x)=(x^2-q^{-4})(x-1)$ since this is the minimal polynomial of $\Rm$ on $\langle v_i\otimes v_j\,|\,i+j\leq 2\rangle$.
Thus if we consider actions of $\,\rho_i=  \sigma_i^2-q^{-4}\,$ and $\,\varepsilon_i=\sigma_i-1\,$ on  $\breve\Vm_{n,2,1}$
we have  $\mathrm{im}(\varepsilon_i)=\mathrm{ker}(\rho_i)$ and $\mathrm{ker}(\varepsilon_i)=\mathrm{im}(\rho_i)$
over $\mathbb Q(q)$ (only if $q^{-4}\neq 1$). The action of  (\ref{def-eqCact}) implies that $e_1\in \mathrm{ker}(\rho_1)
=\mathrm{im}(\varepsilon_1)$. The latter space is spanned by generators $r_1=b_1-\frac {q^2}{[2]}a_{1,2}$, 
$r_2=b_2-\frac {1}{[2]}a_{1,2}$,  as well as $r_j=a_{1,j}-q a_{2,j}$ for $j=3, \ldots, n\,$. Since $e_1$ has to be mapped to 
$\overline b_1$  in the quotient it is thus a linear combinations of the form $e_1=-q^2 r_1+\sum_{i\geq 3}\alpha_ir_i$.

Now the relations in  (\ref{def-eqCact}) for $\lambda=q^{-4}$ also imply that $e_1\in \mathrm{ker}(\varepsilon_i)$ for $i\geq 2$,
which leads to additional constraints that determine the $\alpha_i$ and hence $e_1$ uniquely:
\begin{equation}
 e_1\,=\,- q^2 b_1\,+\,\frac {q^4} {[2]} a_{1,2}\,+\,\frac {q^2} {[2]}\sum_{k\geq 3}^nq^{4-k}(a_{1,k}-qa_{2,k})\;.
\end{equation}
The action of $\sigma_1$  on $\mathbf C_n(q^{-4})$ now implies that
\begin{equation}
 e_2\,=\,q^4\sigma_1.e_1\,=\,\,-q^4 b_2\,+\,\frac {q^4} {[2]} a_{1,2}\,-\,\frac {q^4} {[2]}\sum_{k\geq 3}^nq^{4-k}(a_{1,k}-qa_{2,k}) \;.
\end{equation}
From this it subsequently follows that
\begin{equation}
 \rho_2.e_2\,=\,\frac {q^4-q^{2}} {[2]} \sum_{k\geq 4}^nq^{4-k}(q a_{2,k}+a_{3,k})\;.
\end{equation}
However, by  (\ref{def-eqCact}) we must have $\rho_2.e_2=0$ which leads to a contradiction for $n\geq 4$ and $q^4\neq 1\,$.
\end{proof}

Let us next point out some relations of this lemma to the  the topological construction of
the Temperley-Lieb representation given in Section~6 of \cite{Big03}. 

The identification $q=s$ was motivated
in our case by choosing a fundamental highest weight for quantum-$\sl$ and translates via  (\ref{eq-deftheta})  directly to the
specialization  ${\tt q}{\tt t}=-1$ considered by Bigelow in  \cite{Big03} as well as Lawrence in \cite{Law90}. In terms of 
these variables and pre-composing representations with the involutive automorphism $\invo$ on $B_ n$ given by 
$\invo(\sigma_i)=\sigma_i^{-1}$ we find from Lemma~\ref{lm-TLinLKB} the Temperley-Lieb representation 
$\mathbf{\pi}_{[n-2,2]}$ as the kernel of the following map of $B_n$-modules.
\begin{equation}\label{eq-HCquot}
 \lkb{n}{2}\!\underset{{\tt t}=-{\tt q}^{-1}}{\otimes}\!\mathbb Q({\tt q})\,\longrightarrow\,\mathbf C_n^{\dagger}({\tt q}^{-2})
\;:\quad F_{i,j}\,\mapsto\, e_i\,+\,{\tt q}^{-1}e_j\;.
\end{equation}
Here the action of $B_n$ on $\mathbf C_n^{\dagger}({\tt q}^{-2})$ is given explicitly by 
\begin{equation}\label{def-eqCact-opp}
 \sigma_j.e_j=  e_{j+1}\,,\quad \sigma_j.e_{j+1}= {\tt q}^{2} e_{j}, \quad\mbox{and}\quad  \sigma_j.e_{i}=e_i\;\;\mbox{for } i\not\in\{j,j+1\}\;.
\end{equation}
The Temperley-Lieb representation is found as the kernel of the map (\ref{eq-HCquot}) also by Lawrence (see page 170 in \cite{Law90}), 
however, in the dual or cohomological version of the Lawrence representation. Consequently,  in the homology picture of 
\cite{Law90}  $\mathbf{\pi}_{[n-2,2]}$ is described as a quotient by an $n$-dimensional sub-representation. 

Bigelow finds in Theorem~6.3 of  \cite{Big03} the module $\mathbf{\pi}_{[n-2,2]}$ as the image of $H_2(\tilde Y_{n,2})\otimes R$ in 
 $H_2(\tilde Y_{n,2},\tilde \nu)\otimes R$ by the map induced by the inclusion of pairs, where $\tilde \nu$ is a limit of configurations
in which one of the points of configuration in $\tilde Y_{n,2}$ approaches a puncture or both points approach each other. This
suggests that the module $\mathbf C_n^{\dagger}({\tt q}^{-2})$ is somehow related to the first homology of  $\tilde \nu\,$,
although it is not na\"\i vely obtained from the long exact sequence associated to $(\tilde Y_{n,2},\tilde \nu)$.

The sequence of $B_n$-representations in Lemma~\ref{lm-TLinLKB} fails to split essentially due to the failure of 
$\Imo_1\subset \breve \Vm_1$to split off as a quantum-$\sl$ representation. Again it would be interesting to 
understand this as an obstruction in the context of the topological constructions in \cite{Big03} and  \cite{Law90}
where it contributes to subtle distinctions between various types of homological and cohomological variants of
the LKB-representations. 

More generally, the  $q=s$ specialization of the $\W_{n,k}$ representations will contain the 
$\mathsf{A}_{n,q}\,$-representations $\pi_{[n-k,k]}$ of dimension ${n\choose k}-{n\choose {k-1}}$
as summands by the same arguments used for the case $k=2$ above. This reproduces the 
Temperley-Lieb representations described at the end of Section~5.2 in \cite{Law90}. One may
expect that they are again not direct summands as $B_n$-modules as in the case of $k=2$.

\medskip

The behavior of the representation  $\W_{n,k}$  is very different if we consider them over $\mathbb Q(q,s)$
(where $s-q^{\ell}$ is invertible). In particular, we will show in the following sections that $\W_{n,k}\otimes \mathbb Q(q,s)$
is irreducible for all $n$ and $k\,$. This indicates that the 2-parameter representation over $\mathbb Z[q^{\pm 1},s^{\pm 1}]$
contains significantly more information than the one-parameter specialization discussed above and, especially, the
Temperley-Lieb sub-representation. 

The loss of complexity in the specialization to the Temperley-Lieb representation is exemplified also by the
fact that $\lkb{n}{2}$ is faithful, while the representation $\pi_{[n-2,2]}$ has a non-trivial kernel.  For $\pi_{[2,2]}$
elements in the kernel are specified in Section~3 of \cite{Big02}. 

More complicated elements in the kernel of the Temperley-Lieb representations are constructed in \cite{PT08}. 
In this article Piwocki and Traczyk represent the Temperley-Lieb algebra  $\,\TLsi_n\equiv \mathsf{A}_{n,q}$  
in terms of Kauffman diagrams,   introduce ideals $\Isi_{n,i}$ generated by diagrams with more than $i$ caps and
cups, and consider the kernels of the composite morphism $\Jsi_{n,i}:\,B_n\to \TLsi_n \to \TLsi_{n,i}=\TLsi_n/\Isi_{n,i}\,$. 

In order to relate this to elements in the kernel of $\pi_{[n-2,2]}$ note that
the generator  from  (\ref{eq-Eact}) can be written as $\mathsf E=\mathsf C\circ \mathsf C^{\vee}$ with maps 
\beq
\begin{split}
\mathsf C :\,\mathbb Z[q^{\pm 1}]\to\Imo_1^{\otimes 2} \quad : \quad &
\begin{array}{l}
1\mapsto  v_1\otimes v_0 \,-\, q(v_0\otimes v_1)\,  
\end{array} \\ \rule{0mm}{7mm}
 \mathsf C^{\vee}  :\,\Imo_1^{\otimes 2}\to \mathbb Z[q^{\pm 1}] \quad : \quad &
\begin{array}{l}
\mathsf C^{\vee}(v_1,v_1)=\mathsf C^{\vee}(v_0,v_0)=0, \\
 \mathsf C^{\vee}(v_0,v_1)=1, \mbox{ and } \mathsf C^{\vee}(v_1,v_0)=-q^{-1}.
\end{array} 
\end{split}
\eeq

The action of $\TLsi_n$ on $\Imo_1^{\otimes n} $ can now be extended by associating to planar diagrams with
$a$ start and $b$ end points a map from  $\Imo_1^{\otimes a} $ to $\Imo_1^{\otimes b} $ by assigning 
the tensors $\mathsf C$ and  $\mathsf C^{\vee}$ to cups and caps in respective tensor positions.  Note a diagram
in $\Isi_{n,4}$ must have either at least three cups or three caps. This corresponds to the application, for example,
of three contractions of pairs of tensor factors in  $\Imo_1^{\otimes n} $ with $\mathsf C^{\vee}\,$. Restricted
to  $\,\Imo_1^{\otimes n} \cap \breve \W_{n,2,1}\cong \pi_{[n-2,2]}\,$ all such contractions are zero for degree reasons. 
Similarly, insertion of three or more tensors  with $\mathsf C\,$ cannot have image in $\, \pi_{[n-2,2]}\,$.

We conclude that the ideal $\Isi_{n,4}$ acts trivially on $\, \pi_{[n-2,2]}\,$ and hence, as a braid group representation,
the latter factors through $\,\Jsi_{n,4}:\,B_n\to \TLsi_{n,4}\,$. In  \cite{PT08} Piwocki and Traczyk find a non-trivial
380-crossing 
braid in the kernel of $\Jsi_{9,2}$. Using Theorem~1 in \cite{PT08}  this can be used to construct of a 1520-crossing braid 
$\beta$ in 
the kernel of $\Jsi_{17,4}\,$ and hence also in the kernel of  $\, \pi_{[15,2]}\,$.

Once it is verified that $\beta\neq 1$ (for example, by evaluating it in the LKB representation) this proves that
$\, \pi_{[15,2]}\,$ is not a faithful representation of $B_{17}$. A more accessible candidate may be the 11-crossing braid in $\mathrm{ker}(\Jsi_{11,2})$
which yields a 44-crossing element in $\mathrm{ker}(\Jsi_{21,4})\subseteq\, \mathrm{ker}(\pi_{[19,2]})\, $. We will not engage in the
remaining computations in this article, however, and leave them for future work.

%%%%%%%%%%%%%  End -- Chapter 5   %%%%%%%%%%%%%%%%%%%%
%%%%%%%%%%%%%  Start  -- Chapter 6   %%%%%%%%%%%%%%%%%

\section{Structure of the Verma representations $\Vm_{n,l}$}\label{sec-strVerma}

In this section we look more closely at the structure and decomposition of the Verma module
representations $\Vm_{n,l}$. More specifically, we look 
at eigenspace decompositions of $\Vm_{n,l}$ under the operators
$E^t F^{(t)}$. The main purpose of these decompositions is to allow us to prove
the irreducibility  of the highest weight  representations in the next section.

Recall that we have previously defined $\fRqs = \Q(q, s)$, 
the fraction field of $\Rqs$. In what follows,  we will often speak of $\Vm_{n,l}$ as a vector space over 
$\fRqs$. Of course what we really mean is $\fRqs \otimes_{\Rqs} \Vm_{n,l}$, but we will usually
make no distinction. We could, in the interest of generality, carry out our calculations over a smaller ring, essentially inverting only those elements of $\Rqs$ that
are necessary, but this level of generality adds little to the discussion at hand.

\blm\lb{decomp}
The weight space $\Vm_{n,l}$ splits as 
a $\fRqs [B_n]$-module into a direct
sum of highest weight spaces:
\beq
\Vm_{n,l} = \bigoplus_{k=0}^l F^{(k)} \W_{n,l-k} 
          \cong \bigoplus_{k=0}^l \W_{n,l-k}.
\eeq
\elm
\begin{proof}
We already know that  
$\Vm_{n,l} = \W_{n,l} \oplus \Bm_{n,l} \cong \W_{n,l} \oplus \Vm_{n,l-1}$
as $\Rqs$-modules.
So $\Vm_{n,l}$ does decompose into a direct sum of highest weight spaces
$\bigoplus_{k=0}^l \W_{n,l-k}$. This decomposition does not 
preserve the braid group action, however.

To prove the decomposition
$\Vm_{n,l} = \bigoplus_{k=0}^l F^{(k)} \W_{n,l-k}$ we proceed by induction
on $l$. For $l=0$ we have an obvious identity. Suppose now that
$\Vm_{n,l} = \bigoplus_{k=0}^l F^{(k)} \W_{n,l-k}$ and take
$v = F^{(k)}w \in \Vm_{n,l}$ for some $w \in \W_{n,l-k}$. We apply
$E^t F^{(t)}$ to $v$ to obtain
\begin{align}
E^t F^{(t)} v &= 
         \Bigl[\begin{matrix}t+k\vspace*{-1.6mm} \\k\end{matrix}\Bigr]_q 
                            E^t F^{(t+k)} w \lb{irredeqn} \\
      &= \Bigl[\begin{matrix}t+k\vspace*{-1.6mm} \\k\end{matrix}\Bigr]_q
                                         F^{(k)} \Big( \prod_{j=1}^t 
               (q^{j-k-t} K - q^{k+t-j} K^{-1}) \Big) w \notag \\
      &= \Bigl[\begin{matrix}t+k\vspace*{-1.6mm} \\k\end{matrix}\Bigr]_q 
              \mu_{t,k}^{n,l} \; v \notag
\end{align}
where $\mu_{t,k}^{n,l} \in \Rqs$ is the nonzero constant given by
\beq
\mu_{t,k}^{n,l}= \prod_{j=1}^t (s^n q^{-2l+k-t+j} - s^{-n} q^{2l-k+t-j}).
\eeq
Thus, in particular $EF^{(1)}v = [k+1]_q \mu_{1,k}^{n,l} v$. Since
the constants $[k+1]_q \mu_{1,k}^{n,l} \in \Rqs$ are
distinct for distinct $k$, we see that
the decomposition 
$\Vm_{n,l} = \bigoplus_{k=0}^l F^{(k)} \W_{n,l-k}$
is the eigenspace decomposition of the transformation
$EF^{(1)}$. The eigenvalues $[k+1]_q \mu_{1,k}^{n,l}$ are 
each nonzero, so we see that
the map $F^{(1)} : \Vm_{n,l} \to \Vm_{n,l+1}$ is injective. The image
of this map (over the fraction field) is 
$\mbox{Im}(F^{(1)}) = \bigoplus_{k=1}^{l+1} F^{(k)}\W_{n,l+1-k}$
and it is clear from \eqref{irredeqn}
that $\mbox{Im}(F^{(1)}) \cap \W_{n,l+1} = 0$. Counting dimensions,
we see that $\Vm_{n,l+1} = \bigoplus_{k=0}^{l+1} F^{(k)} \W_{n,l+1-k}$.
\end{proof}

Having obtained a decomposition of
$\Vm_{n,l}$, we would now like to obtain
a similar decomposition of the highest weight spaces
$\W_{n,l}$ by restricting the braid action. 

Consider the $B_{n+1}$-action on $\W_{n+1,l}$. The
map $\Vm^{\otimes n} \to \Vm^{\otimes (n+1)}$ defined by
$v_{\vec \alpha} \mapsto v_0 \otimes  v_{\vec \alpha}$ 
gives us an inclusion
$\W_{n,l} \into \W_{n+1,l}$. In the standard basis of $\W_{n+1,l}$
the elements of $\W_{n,l}$ correspond to the vectors
$\Phi(a_{\vec{\alpha}})$ where
$\vec \alpha = (\alpha_j, \dots, \alpha_n)$ for $j > 2$ (see \eqref{Abasis} and
\eqref{map}). 
We also have the inclusion
$B_{n} \into B_{n+1}$ that takes $\sigma_i \in B_{n}$ to
$\sigma_{i+1} \in B_{n+1}$. With this identification the inclusion
$\W_{n,l} \into \W_{n+1,l}$ is $B_{n}$-equivariant. The quotient
$\W_{n+1,l} / \W_{n,l}$ is isomorphic to $\Vm_{n,l-1}$ as
an $\Rqs [B_{n}]$-module. The isomorphism is given by
\beq\lb{quotient map}
\Phi(a_{\vec{\alpha}}) \mapsto v_{\vec{\alpha}}.
\eeq
Let $\psi : \W_{n+1,l} \to \Vm_{n,l-1}$ be the composition of the 
quotient map $\W_{n+1,l} \to \W_{n+1,l} / \W_{n,l}$ with the isomorphism
given in \eqref{quotient map}. We seek a splitting of $\psi$. 
\bdf\lb{def-alphak}
Let $c_{k,j} \in \Rqs$ be recursively defined by setting
$c_{k,0} = 1$ and
\beq
 c_{k,j+1}
=  \frac{s^{-n-1}q^{2l-k+j-1} - s^{n+1}q^{-2l+k-j+1}}{s^{n} q^{-2(l-k)}}
c_{k,j}
\eeq
For each $k = 1,2, \dots,l$ we define a map 
$\alpha_k : \W_{n,l-k} \to \Vm_{n+1,l}$ by
\beq
\alpha_k : w \mapsto \sum_{j=0}^k c_{k,j} F^{(k-j)}( v_j \otimes w).
\eeq
\edf
Let us take $w \in \W_{n,l-k}$ 
and compute the action of $E$ on $\alpha_k (w)$:
\begin{align*}
E \alpha_k (w) &= \sum_{j=0}^k c_{k,j} E F^{(k-j)}( v_j \otimes w)\\
               &= \sum_{j=0}^k c_{k,j} \bigl( F^{(k-j)} E 
       + F^{(k-j-1)} (q^{1-k+j}K - q^{k-j-1}K^{-1}) \bigr) 
                                                 ( v_j \otimes w) \\
               &= \sum_{j=1}^k c_{k,j}
           s^{n}q^{-2(l-k)} F^{(k-j)}(v_{j-1} \otimes w) \\
               & \phantom{adad} + \sum_{j=0}^{k-1} c_{k,j} 
     \bigl( s^{n+1} q^{-2l +k-j+1} - s^{-n-1} q^{2l-k+j-1} \bigr)
                     F^{(k-j-1)}(v_j \otimes w)  \\
               &= \sum_{j=0}^{k-1} 
         \bigl( c_{k,j+1} s^{n}q^{-2(l-k)} 
         + c_{k,j} ( s^{n+1} q^{-2l +k-j+1} - s^{-n-1} q^{2l-k+j-1}) \bigr)
                    F^{(k-j-1)}(v_j \otimes w)  \\
               &= 0.
\end{align*}
Thus, $\alpha_k$ actually maps $\W_{n, l-k}$ into 
$\W_{n+1,l} \subset \Vm_{n+1,l}$. Notice, in the last equality we see the reason behind the 
definition of the coefficients $c_{k,j}$ in Definition \ref{def-alphak}. 
Namely, they have been  defined to allow $E \circ \alpha_k$ to vanish on $\W_{n,l-k}$.

In the standard basis of
$\W_{n+1,l}$  the element $\alpha_k w$ corresponds, modulo $\W_{n,l}$,
to a multiple of $v_1 \otimes F^{(k-1)}w$. To be
more precise:
\beq\lb{alpha}
\alpha_k w = \lambda_k \Phi( v_1 \otimes F^{(k-1)}w ) 
                              \mbox{   mod } \W_{n,l} 
\eeq
where
\begin{align}
\lambda_k &= s^{1-k} q^{k-1} (s-s^{-1}) \nonumber
                     + c_{k,1} s^{1-k} q^{2k-2} \\
   &= s^{-2n-k} q^{4l-k-3} - s^{-k} q^{k-1}.
\end{align}
Thus, using the identification 
$\Vm_{n,l-1} = \bigoplus_{k=0}^{l-1} F^{(k)} \W_{n,l-1-k}
          \cong \bigoplus_{k=1}^l \W_{n,l-k}$
given by
Lemma \ref{decomp} we see that 
$\psi \circ \alpha_k$ acts on $\W_{n, l-k}$ as multiplication
by the nonzero constant $\lambda_k$.
\bdf
Define a map $\alpha : \Vm_{n,l-1} \to \W_{n+1,l}$ by 
\beq
\alpha = \bigoplus_{k=1}^{l} \lambda_k^{-1} \alpha_k.
\eeq
\edf
The previous discussion yields the following:
\blm\lb{splitting}
The map $\alpha$ defines a $B_{n}$-equivariant splitting
of the map $\psi : \W_{n+1,l} \to \Vm_{n,l-1}$. This gives
a decomposition as $B_{n}$-modules
\beq\lb{decomposition}
\W_{n+1,l} = \bigoplus_{k=0}^l \W_{n,l-k}.
\eeq
\elm

%%%%%%%%%%%%%  End -- Chapter 6   %%%%%%%%%%%%%%%%%
%%%%%%%%%%%%%  Start -- Chapter 7   %%%%%%%%%%%%%%%%%

\section{Irreducibility of the Representations}\label{sec-Irred}

In this last section we wish to prove Theorem~\ref{thm-irred}, namely
that the highest weight representations $\W_{n,l}$ are irreducible
over the fraction field $\fRqs$. The proof
makes use of the decompositions of the previous section
and proceeds by induction on $n$.
Notice that, in the general case, if
$C \subset \W_{n,l}$ is a $B_n$-submodule, then as a
$B_{n-1}$-module it must decompose into 
a direct sum of lower degree submodules following the decomposition
\[
\W_{n,l} = \bigoplus_{j=0}^l  \W_{n-1,j}\;.
\]
By the induction hypothesis, each of these summands is an irreducible
representation of $B_{n-1}$ so that $C$ must be a direct 
sum of some collection of these
$\W_{n-1,j}$ (for more detail see the proof at the end of the section). 
In what follows we give explicit computations of 
the action of 
$\sigma_1 \in B_n$ on certain elements of these
components. These computations show that 
we must, in fact, have $\W_{n-1, j} \subset C$ for all $j$,
thus proving the theorem.

To start, let us suppose  that $v \in \Vm_{n,l}$, then by 
Lemma~\ref{decomp}
we have $v = w_0 + F^{(1)} w_1 + \cdots + F^{(l)} w_l$
for some $w_t \in \W_{n,l-t}$. We would like to be able to 
describe these vectors $w_t$ in terms of $v$. 

For any $t \leq l$ we apply $E^t$ to $v$
to obtain
\begin{align}
E^t v &=  E^t F^{(t)} w_t + E^t F^{(t+1)} w_{t+1} 
                       + \cdots + E^t F^{(l)} w_l \\
  &= \mu_{t,0}^{n,l-t} w_t + \mu_{t,1}^{n,l-t} F^{(1)} w_{t+1} 
                       + \cdots + \mu_{t, l-t}^{n,l-t} F^{(l-t)} w_l
\end{align}
so that we can solve recursively for $w_t$:
\beq\lb{eq-mu}
w_t = \frac{1}{\mu_{t,0}^{n,l-t}} 
   \Big( E^t v - \mu_{t,1}^{n,l-t} F^{(1)} w_{t+1}
                       - \cdots - \mu_{t, l-t}^{n,l-t} 
                            F^{(l-t)} w_l \Big).
\eeq
Proceeding by induction, we see that we must have
\beq\lb{eq-w}
w_t = \sum_{i=0}^{l-t} z_{t,i}^{n,l} F^{(i)} E^{t+i} v
\eeq
for some coefficients $z_{t,i}^{n,l} \in \fRqs$. 
We see from \eqref{eq-mu} that
$z_{t,0}^{n,l} = 1/ \mu_{t,0}^{n,1-t}$ and 
an induction argument shows that, in general,
\beq\lb{eqn-z}
z_{t,i}^{n,l}|_{q=1} = (-1)^i z_{t+i,0}^{n,l}|_{q=1} 
= (-1)^i (s^n-s^{-n})^{-t-i}.
\eeq
In particular, the coefficients $z_{t,i}^{n,l}$ are never zero.

\bex\label{example}
Let us define $\nu_j = v_j \otimes v_0^{\otimes (n-1)} \in \Vm_{n,j}$. 
Then, as above,
\beq\lb{eqn-vj}
\nu_j= w_{j,0} + F^{(1)} w_{j,1} + \dots + F^{(j)} w_{j,j}
\eeq 
for 
some $w_{j,i} \in \W_{n,j-i}$. Let us use \eqref{eq-w} to define
\beq\lb{eqn-wj}
\omega_j \stackrel{\text{def}}{=} w_{j,0} 
     = \sum_{i=0}^{j} z_{0,i}^{n,j} F^{(i)} E^i \nu_j.
\eeq 
In other words, $\omega_j$ is the first term of $\nu_j$ 
in the decomposition $\Vm_{n,j} = \bigoplus_{k=0}^j \W_{n,j-k}$.
Since $E^i \nu_j = s^{(n-1)i} \nu_{j-i}$, we see from \eqref{eqn-wj}
that $\omega_{j} \neq 0$ for all $j$. Also, from \eqref{eqn-vj} we see 
that $w_{j,i} = s^{(n-1)i} z_{i,0}^{n,j} \omega_{j-i}$. Thus, equation
\eqref{eqn-vj} can be written as
\beq
\nu_j = \sum_{i=0}^j s^{(n-1)i} z_{i,0}^{n,j} \, F^{(i)} \omega_{j-i}.
\eeq
\eex

\blm\lb{lem-irr}
Let us define $\nu_{j,k} = v_j \otimes F^{(k)} v_0^{\otimes (n-1)} \in \Vm_{n,j+k}$
with $j+k \leq l$.
Then
\[
\nu_{j,k} = \sum_{i=0}^{j+k} \Gamma_{j,k,i} F^{(i)} \omega_{j+k-i}
\]
where $\Gamma_{j,k,i} \in \fRqs$ such that 
$(1-s^{2n})^l \Gamma_{j,k,i} |_{q=1}$ is a 
Laurent polynomial in $s$ with smallest degree term given by
\beq\lb{eqn-Gamma}
\begin{cases}
\binom{j+k-i}{k-i} s^i & 0 \leq i \leq k \\
(-1)^{k-i} \binom{i}{i-k} s^{i + 2(i-k)(n-1)} & k < i \leq  j+k.
\end{cases}
\eeq
\elm
\begin{proof}
From the previous example we have $\nu_j = \nu_{j,0}$ so that
$\Gamma_{j,0,i} = s^{(n-1)i} z_{i,0}^{n,j}$ and it is easy to verify
the lemma for the case $k=0$.

In the general case, we first notice that $\nu_{j,k}$ can be expressed 
as follows:
\beq\lb{eqn-nu}
\nu_{j,k} = \sum_{r=0}^k \gamma_{j,k,r} F^{(k-r)} \nu_{j+r}
\eeq
where the coefficients $\gamma_{j,k,r}$ are defined recursively by 
first setting
\beq
\gamma_{j,0,r} =
\begin{cases}
1, &\text{if $r=0$;}\\
0, &\text{if $r\neq 0$}
\end{cases}
\eeq
and then defining for $k \geq 1$
\beq\lb{e-gamma}
\gamma_{j,k,r} = \frac{q^{-2j}}{[k]_q} \Big(
 [k-r]_q s \gamma_{j,k-1,r} - [j+1]_q (s^2 q^{-j} - q^j ) \gamma_{j+1,k-1,r-1}
\Big).
\eeq
The verification of this fact follows by an induction argument
from the identity
\beq
\nu_{j,k} = \frac{sq^{-2j}}{[k]_q} F^{(1)} \nu_{j,k-1}
  -\frac{q^{-2j}[j+1]_q}{[k]_q}(s^2 q^{-j} - q^j) \nu_{j+1,k-1}.
\eeq
Using the $k=0$ case, equation \eqref{eqn-nu} becomes
\begin{align}
\nu_{j,k} &= \sum_{r=0}^k \sum_{t=0}^{j+r}
            \gamma_{j,k,r} \Gamma_{j+r,0,t} 
        \Bigl[\begin{matrix}k-r+t\vspace*{-1.6mm} \\t\end{matrix}\Bigr]_q
              F^{(k-r+t)} \omega_{j+r-t} \\
         &= \sum_{i=0}^{j+k} \Gamma_{j,k,i} F^{(i)} \omega_{j+k-i}
\end{align}
where
\[
\Gamma_{j,k,i} \stackrel{\text{def}}{=} 
    \sum_{ \{r,t \mid k-r+t = i \} } \gamma_{j,k,r} \Gamma_{j+r,0,t}
        \Bigl[\begin{matrix}i\vspace*{-1.6mm} \\t\end{matrix}\Bigr]_q.
\]

Now, we would like to know something about these coefficients
$\Gamma_{j,k,i}$. At least we would like to know that they are nonzero.
The relation in \eqref{e-gamma} along with an easy induction argument
show that $\gamma_{j,k,r} |_{q=1}$ is a polynomial in $s$ of the form
\[
\gamma_{j,k,r} |_{q=1} = \binom{j+r}{r} s^{k-r} + (\text{higher degree terms}).
\]
Also, from \eqref{eqn-z} we have
\[
\Gamma_{j+r,0,t} |_{q=1} = s^{(n-1)t}(s^n - s^{-n})^{-t}.
\]
So for each $r,t$ with $k-r+t=i$ we see that setting $q=1$ in
$(1-s^{2n})^l \gamma_{j,k,r} \Gamma_{j+r,0,t}
        \Bigl[\begin{matrix}i\vspace*{-1.6mm} \\t\end{matrix}\Bigr]_q$
will indeed give us a Laurent polynomial in $s$ with smallest degree term
\[
(-1)^t \binom{i}{t} \binom{j+r}{r} s^{i + 2t(n-1)}.
\]
Since the degree of this term is positively related to $t$, the
overall smallest degree term of $\Gamma_{j,k,i}$ will occur when 
$t$ is as small as possible. For $0 \leq i \leq k$, the smallest $t$
may be is 0 and in this case we also have $r = k-i$. For $k < i \leq j+k$
the smallest $t$ may be is $i-k$ and in this case we have $r=0$.
This proves \eqref{eqn-Gamma} and the lemma.
\end{proof}

Suppose we take $w \in \W_{n,l-k}$ to be the basis vector given by
\begin{align}
w &= \Phi (v_1 \otimes u) \\
   &= \sum_{t=0}^{l-k} b_t \, v_t \otimes E^{t-1} u
\end{align}
for some $u \in \Vm_{n-1,l-k-1}$ and where the 
coefficients $b_t$ are given as in \eqref{map}. We then have
\[
\alpha_k w = \sum_{j=0}^k \sum_{t=0}^{l-k} c_{k,j} b_t
              F^{(k-j)} ( v_j \otimes v_t \otimes E^{t-1} u ).
\] 
Let us (temporarily) set 
\[
d_{h,j,t} = q^{h(h-1)/2} s^{-(j+t)} q^{2(j-h)(t+h)}. 
\]
We act on $\alpha_k w$ by $\sigma_1 \in B_{n+1}$ and compute
\beq\lb{e-example}
\sigma_1 (\alpha_k w) = \sum_{h=0}^{\infty} \sum_{j=0}^k \sum_{t=0}^{l-k}
   c_{k,j} \, b_t \, d_{h,j,t} \,
F^{(k-j)} \Big( F^{(h)} v_t \otimes E^h v_j \otimes E^{t-1}u \Big).
\eeq
Recall that $\psi : \W_{n+1,l} \to \Vm_{n,l-1}$ is the map
that first mods out by $\W_{n,l}$, then projects to $\Am_{n+1,l}$, then
removes the leading $v_1$ component in the tensor product. 
So applying $\psi$ to $\sigma_1 (\alpha_k w)$, the only terms to survive
are those for which $h+t \leq 1$, and we obtain
\beq\lb{e-braidcalc}
\psi(\sigma_1 (\alpha_k w)) 
   = \sum_{j=0}^k 
  F^{(k-j)}(\eta_j \, v_j \otimes u + b_0 \kappa_j \, v_{j-1} \otimes E^{-1}u).
\eeq
where the coefficients $\eta_j$ and $\kappa_j$ are calculated to be
\begin{align}
\eta_j &=c_{k,j} s^{-(k+1)} q^{2k}, \\
\kappa_j &= s^{-k} q^{2(k-1)} (s-s^{-1}) 
               \Big( c_{k,j}  + c_{k,j-1} s q^{2-j-k} \Big)\,.
\end{align}

Hence we now have formulae for the $\sigma_1$-action on
$\W_{n,l}$ in terms of the decomposition \eqref{decomposition}. 
We make use of this in the next lemma which will be our main tool
in proving Theorem \ref{thm-irred}.

\blm\lb{lm-irrmain}
Let $0 < k < l$ and consider $\omega_{l-k} \in \W_{n,l-k}$
as given in \eqref{eqn-wj}.
Then  in terms of the
decomposition \eqref{decomposition},
$\sigma_1 (\alpha_k \omega_{l-k})$ 
has nontrivial components in $\W_{n,l-r}$ for
all $r = 1, 2, \dots, k+1$.
\elm
\begin{proof}
Let us take
$u = F^{(l-k-1)}(v_0^{\otimes (n-1)}) \in \W_{n-1,l-k-1}$ and define
\[
\mathbf{w} = \sum_{t=1}^{l-k} 
        b_t  v_t \otimes E^{t-1} u + x_0 v_0 \otimes F^{(1)}u
\]
where
\[
x_0 = \frac{b_0}{[l-k]_q \mu_{1,l-k-1}^{n-1,l-k-1}}.
\] 
Then $\mathbf{w} \in \W_{n,l-k}$, which follows by the
computations found in \eqref{irredeqn} and the proof of 
Theorem \ref{thm-freemodule},
and comparing the expressions of $\mathbf{w}$ and $\omega_{l-k}$
in the standard basis we see that $\omega_{l-k}$ is a nonzero
multiple of $\mathbf{w}$. Hence if 
we prove the lemma for $\mathbf{w}$ then
it will also follow for $\omega_{l-k}$.

Notice that
the discussion following Lemma \ref{splitting} will apply 
formally to $\mathbf{w}$
if we make the substitutions $b_0 \to x_0$ and $E^{-1} \to F^{(1)}$.
Thus in the present case \eqref{e-braidcalc} becomes
\beq\lb{e-braidcalc2}
\psi(\sigma_1 (\alpha_k \mathbf{w}))
   = \sum_{j=0}^k \eta_j \, F^{(k-j)}(v_j \otimes u)
     + \sum_{j=1}^k x_0 \kappa_j \, F^{(k-j)}(v_{j-1} \otimes F^{(1)} u).
\eeq
We have
\begin{align}
F^{(k-j)} (v_{j-1} \otimes F^{(1)}u) &= 
             s q^{-2(j-1)} [k-j+1]_q F^{(k-j+1)}(v_{j-1} \otimes u) \nonumber \\
 & \phantom{adad} - s q^{-2(j-1)} 
      (s q^{1-j} s^{-1} q^{j-1}) [j]_q F^{(k-j)}(v_j \otimes u)
\end{align}
which allows us to write \eqref{e-braidcalc2} as
\beq
\psi(\sigma_1 (\alpha_k \mathbf{w})) = 
        \sum_{j=0}^k \Upsilon_j F^{(k-j)}(v_j \otimes u)
\eeq
where
\beq
\Upsilon_j = \eta_j + x_0 s q^{-2(j-1)} \big(
           \kappa_{j+1} q^{-2} [k]_q  
           - \kappa_j (s q^{1-j} - s^{-1} q^{j-1}) [j]_q \big) 
\eeq
which makes sense for all $j= 0,1,\dots,k$ so long as we 
define $\kappa_0 = \kappa_{k+1} = 0$.

We apply Lemma \ref{lem-irr} to obtain
\beq
\psi(\sigma_1 (\alpha_k \mathbf{w})) 
   = \sum_{j=0}^k \sum_{i=0}^{j+l-k-1} \Upsilon_j \Gamma_{j,l-k-1,i}
    \Bigl[\begin{matrix}k-j+i\vspace*{-1.6mm} \\i\end{matrix}\Bigr]_q
 F^{(k-j+i)} \omega_{j+l-k-1-i}.
\eeq
Thus, if $\psi(\sigma_1 (\alpha_k \mathbf{w})) = 
\sum_{r=0}^{l-1} F^{(r)} \mathbf{w}_r$
where $\mathbf{w}_r \in \W_{n,l-1-r}$, then for
all $r=0,1,\dots, k$ we will have
\beq\lb{eqn-irrcoeff2}
\mathbf{w}_r = \sum_{i=0}^r \Upsilon_{k-r+i} \Gamma_{k-r+i,l-k-1,i}
\Bigl[\begin{matrix}r\vspace*{-1.6mm} \\i\end{matrix}\Bigr]_q
\omega_{l-1-r}.
\eeq
Thus, to complete the proof we need only show that the coefficient in
\eqref{eqn-irrcoeff2} is nonzero for all $r=0,1,\dots,k$. 
But we can check that
$s^{n+2(k+1)}(s^{n-1}-s^{-n+1}) \Upsilon_{k-r+i} |_{q=1}$ is
a Laurent polynomial in $s$ having smallest degree term
\[
\frac{k}{l-k} s^{-(i+k-r)(2n+1)}.
\]
Compare this to the minimum degree term of 
$\Gamma_{k-r-i, l-k-1,i}$
given via \eqref{eqn-Gamma} and we see that, after multiplying by an
appropriate constant not depending on $i$, the smallest degree term of 
$\Upsilon_{k-r+i} \Gamma_{k-r+i,l-k-1,i}$ is a strictly decreasing
function of $i$ (details are left to reader). Thus, 
$\sum_{i=0}^r \Upsilon_{k-r+i} \Gamma_{k-r+i,l-k-1,i}
\Bigl[\begin{matrix}r\vspace*{-1.6mm} \\i\end{matrix}\Bigr]_q$ must be 
nonzero.
\end{proof}

Before we prove the irreducibility of $\W_{n,l}$ we need one 
last result.

\blm\lb{wmax}
For any $w \in \W_{n,l}$ there is a polynomial 
$P_w (x) \in \Rqs [x]$
such that $P_w (\sigma_1) w \in \Rqs \cdot w_{max}$ where
\[
w_{max} = \Phi ( v_1 \otimes v_{l-1} \otimes v_0^{\otimes (n-2)}).
\] 
\elm
\begin{proof}
We give the standard basis of $\W_{n,l}$ the following ``lexicographical''
ordering:
\[
w_{\vec \alpha} < w_{\vec \gamma} \mbox{   if   }
\begin{cases}
|\vec \alpha| < |\vec \gamma| \mbox{  or}\\
|\vec \alpha| = |\vec \gamma| \mbox{ and } \vec \alpha < \vec \gamma
\end{cases}
\] 
where by
$| \vec \alpha |$ we mean the number of components in the multi-index (that is, if
$\vec \alpha = (\alpha_j, \dots, \alpha_n )$ then $| \vec \alpha | = n-j+1$)
and by $\vec \alpha < \vec \gamma$ we mean the usual lexicographical 
ordering on ordered tuples of integers that have the same number of components.  

The element $w_{max}$ is the maximal 
element in this ordering and the braid $\sigma_1$ acts on $w_{max}$
as the invertible constant 
\[
 (-1)^l s^{-2l} q^{l(l-1)}.
\]
Now take $w_{\vec \alpha} < w_{max}$. A simple case-by-case analysis
shows that there is a polynomial $P_{\vec \alpha} (x)$ of degree at most 2
such that $P_{\vec \alpha} (\sigma_1) w_{\vec \alpha}$ either belongs to
$\Rqs \cdot w_{max}$ or is a sum of strictly higher order terms.

For instance, suppose $|\vec \alpha | = n-1$ and consider 
\[
w = \Phi ( v_1 \otimes v_k \otimes u )
\]
for some $k > 0$ and some $u \in \Vm_{n-2, l-k-1}$. Then 
we compute
\[
\sigma_1 w = 
     \sum_{i \geq 0} \Phi (v_1 \otimes v_{k+i} \otimes u_i)
\] 
for some $u_i \in \Vm_{n-2, l-k-i-1}$. When $i=0$ we have
$u_0 = z_0 u$ for a nonzero constant $z_0 \in \Rqs$
so that we will have 
\[
(\sigma_1 - z_0) w =
     \sum_{i \geq 1} \Phi (v_1 \otimes v_{k+i} \otimes u_i).
\]
Similar arguments apply in the remaining cases. 

Thus, we can multiply any basis vector by a polynomial in $\sigma_1$
and get a sum of higher order terms. The result now follows by induction and
by commutativity of polynomials in $\sigma_1$. 
\end{proof}

We now come to the main result of this section.

\btm
The $B_n$ representations $\W_{n,l}$ are irreducible
over the fraction field  $\fRqs$.
\etm
\begin{proof}
We proceed by induction on $n$. The base case when $n=2$ is trivial
since the dimension of $\W_{2,l}$ is 1 for all $l \geq 0$. Suppose
now the theorem is true for all $k < n$. 
Suppose $C \subset \W_{n,l}$ is a $B_n$-submodule. As a 
$B_{n-1}$-module we have seen that $\W_{n,l}$ decomposes into
a direct sum
\[
\W_{n,l} = \bigoplus_{j=0}^l  \W_{n-1,j}\;.
\]
By the induction hypothesis, each of these summands is an irreducible
representation of $B_{n-1}$ and since the dimensions are different they are
inequivalent. (Note, in case n=3 the summands are isomorphic, but the braid
action is given by distinct eigenvalues on distinct summands.)
Thus $C$ must be a direct sum of some collection of these
$\W_{n-1,j}$.

It is clear that for any element $w \in \W_{n-1,l}$ there
is a braid $\beta \in B_n$ such that $\beta w$ represents
a nonzero class in $\W_{n,l} / \W_{n-1,l}$. 
%
%
%
\begin{comment}
Hence, we have
\[
\alpha \circ \psi (\beta w) = \beta w \mbox{   mod } \W_{n-1,l}
\]
which means that 
$\mbox{Image}(\alpha) \cap \fRqs [B_n ] \cdot \W_{n-1,l} \neq 0$.
\end{comment}
Hence, we can assume $\W_{n-1,j} \subset C$ for some $j < l$.
But Lemma \ref{lm-irrmain} shows that we must, in fact, have
$\W_{n-1,j} \subset C$ for all $j < l$.

Hence, to complete the proof all that remains is to show 
that $\W_{n-1,l} \subset C$.
But by Lemma \ref{wmax} 
any $w \in C$ can be multiplied by a polynomial $P_w (\sigma_1)$
to obtain a multiple of $w_{max}$. For instance,  if we take
$w$ to be the image under $\alpha$ of the
maximal basis element in  $\W_{n-1, l-1}$, then
$P_w (\sigma_1) w$ will be a nonzero multiple of $w_{max}$. But
$\sigma_1^{-1} \sigma_2^{-1} (w_{max}) \in \W_{n-1,l}$ so we see that
$C$ must equal $\W_{n,l}$.
\end{proof}

%%%%%%%%%%%%%  End -- Chapter 6   %%%%%%%%%%%%%%%%%

\end{document}